\documentclass{article}
\usepackage{amsmath,amsthm,amsfonts,latexsym}

\begin{document}

\newtheorem*{theo}{Theorem}
\newtheorem*{pro} {Proposition}
\newtheorem*{cor} {Corollary}
\newtheorem*{lem} {Lemma}
\newtheorem{theorem}{Theorem}[section]
\newtheorem{corollary}[theorem]{Corollary}
\newtheorem{lemma}[theorem]{Lemma}
\newtheorem{proposition}[theorem]{Proposition}
\newtheorem{conjecture}[theorem]{Conjecture}

 \newtheorem{definition}[theorem]{Definition}
 \newtheorem{remark}[theorem]{Remark}

\newcommand{\Naturali}{{\mathbb{N}}}
\newcommand{\Reali}{{\mathbb{R}}}
\newcommand{\Complessi}{{\mathbb{C}}}
\newcommand{\Toro}{{\mathbb{T}}}
\newcommand{\Relativi}{{\mathbb{Z}}}
\newcommand{\HH}{\mathfrak H}
\newcommand{\KK}{\mathfrak K}
\newcommand{\LL}{\mathfrak L}
\def\A{{\cal A}}
\def\B{{\cal B}}
\def\C{{\cal C}}
\def\H{{\cal H}}
\def\L{{\cal L}}
\def\N{{\cal N}}
\def\M{{\cal M}}
\def\O{{\cal O}}
\def\R{{\cal R}}

\def\V{{\cal V}}
\def\rop{{\rm op}}
\def\D{\Delta}
\def\U{U^{\alpha}}
\def\u{u^{\alpha}}
\def\f{\varphi}
\def\e{\varepsilon}
\def\a{\alpha}
\def\tR{\tilde R}
\def\tS{\tilde S}

\title{On amenability and  co-amenability of algebraic quantum groups
and  their corepresentations}

\author{Erik B\'edos$^*$ \ \ \
Roberto Conti \ \ \
Lars Tuset$^*$}
\date{October 2001 - with some minor changes in November 2001}

\maketitle

\markboth{R. Conti, E. B\'edos, L. Tuset}{Amenability and algebraic
quantum groups}
\renewcommand{\sectionmark}[1]{}

\begin{abstract}
We introduce and study  several amenability properties for
unitary corepresentations and $*$-representations of algebraic quantum
groups,
which may be used to characterize amenability or co-amenability of
such groups. As a background for this study, we also investigate the
involved tensor C$^{*}$-categories.

\vskip 1.5cm

\noindent \emph{Subj. Class.}: Quantum groups,
representations and corepresentations,
amenability.



\end{abstract}

\vfill
\thanks{\noindent $^*$  Partially supported by the Norwegian Research
Council.}

\newpage

\section{Introduction}
The concept of amenability plays an important role in the theory of
locally compact groups and in the theory of operator algebras (see
\cite{Pa} and references therein). Quite naturally,
the concept of
amenability and its companion, co-amenability, have been introduced
and studied by several authors in
various settings
related to quantum groups (see \cite{Vo, ES, Ruan, BS, Ba, Ba2, Bl,
Ng1, Ng2, BMT1, BMT2, BMT3}, in chronological order).

In this paper, we introduce several
concepts of ``amenability'' for unitary corepresentations of analytic
extensions of
algebraic quantum groups (as defined by J.\ Kustermans and A.\ van Daele
\cite{KuVD}) :
\emph{co-amenability} (a notion inspired by results in \cite{BMT1,
BMT2}), \emph{amenability} (inspired by the concept of amenability of a
unitary representation of a locally compact group introduced by M. Bekka
\cite{Bek}) and  the \emph{weak containment property} (inspired by the
classical
characterization of the amenability of a group in terms of weak
containment).

We present  several equivalent formulations of these properties,
and use them  to  characterize amenability or co-amenability of  algebraic
quantum
groups.
After some preliminaries in Section 2, we begin the paper with a categorical
interlude in Section 3 where we show how the
category of non-degenerate $*$-representations of the universal
C*-algebraic quantum group associated to an algebraic quantum
group \cite{Kus} and the category of unitary corepresentations of the
analytic
extension of the dual quantum group (with opposite co-product) are
naturally isomorphic as tensor C*-categories. This result enables to
transfer all notions introduced for unitary corepresentations to
non-degenerate $*$-representations and vice-versa. We also derive the
absorbing property of the fundamental multiplicative unitary and of
the regular representation. In Section 4, we
introduce the conjugate corepresentation and the Hilbert-Schmidt
corepresentation associated with a unitary corepresentation. Section 5
is devoted to co-amenability, Section 6 to amenability and  Section 7 to
the weak containement property, where we also consider briefly property (T)
for algebraic quantum groups. In Section 8, we gather some remarks on
the relationship between these amenability concepts. Finally, in Section 9,
we
specify our study of amenability to the setting of algebraic quantum groups
of discrete
type, where it is possible to exploit the structural properties of
these quantum groups to push our analysis further.

\vspace{2ex}

\noindent
Every vector space
will be over the ground field $\mathbb C$.
Given a set $V$, $\iota_V$ denotes the identity map on it
(but we simply write $\iota$ when there is no danger of confusion).
If $\H$ is a Hilbert space, then
$B(\H)$ (resp. $ B_0(\H)$)
denotes the algebra of all bounded
(resp. compact)
linear operators acting on $\H.$
If $\B$ is a $*$-algebra,
$M(\B)$ denotes the multiplier algebra of $\B$.
If $\B$ is unital, we denote its unit by $I_\B$, or by $I$
when this causes no confusion.
In this case,
we set $\cal{U}(\B)$ for the unitary group of $\B$.
We denote by $S(\B)$ the state space of the $C^*$-algebra $\B$.
As usual $\otimes$ denotes tensor product;
depending on the context,
it may be
the tensor product of vector spaces,
the Hilbert space tensor product
or the minimal (that is, spatial) tensor product of $C^*$-algebras,
$\bar \otimes$ being used for tensor products in the von Neumann
algebra setting.
However, we often use $\odot$ to stress that we are dealing with
an algebraic tensor product.
If $V, W$ are vector spaces,
$\chi: V \otimes W \to W \otimes V$
is the flip map sending $v \otimes w$ to $w \otimes v$
($v \in V, w \in W$);
if $\H$ is a Hilbert space
then $\Sigma$ is the flip map on $\H \otimes \H$.
We use the leg-numbering notation as introduced in \cite{BS}.

\section{Preliminaries}

Throughout this paper, $(\A,\Delta)$
 denotes an algebraic quantum group in the sense of \cite{VD1},
 see also \cite{VD2,KuVD},
 where $\Delta: \A \to M(\A \odot \A)$ is the co-product map.
 We follow notation and use terminology from these papers.
Hence,
$S$ denotes the antipode of $(\A,\Delta),$ $\varepsilon$ its counit and
$\varphi$ is a fixed faithful left Haar functional. This  functional $\f$
s not necessarily tracial (or central).
However, there is a unique bijective homomorphism $\rho \colon \A \to \A$
such that $\f(a b)=\f(b \rho(a))$, for all $a,b\in \A$. Moreover,
$\rho(\rho(a^*)^*)=a$.

\vspace{1ex}
The pair $(\A_{r},\Delta_{r})$
denotes the associated analytic extension
(which is a reduced locally compact quantum group in the sense of
\cite{KV}), $\pi_r: \A \to \A_r \subset  B(\H)$
 is the (left) regular representation of $\A$
acting on the GNS Hilbert space $\H$ of $\varphi$, $\Lambda : \A \to
\H$ is the canonical injection,
$W \in M(\A_{r} \otimes B_{0}(\H))$ is the
associated multiplicative unitary,
$\M = \A_r '' = \pi_r(\A)''$ is the von Neumann algebra generated by
$\pi_r(\A)$ and  $R$ is the anti-unitary antipode
(which is defined on $\M$).

\vspace{1ex}

We denote by $(\hat \A,\hat \Delta)$ the dual algebraic quantum group
and by $(\hat \A_r,\hat \Delta_{r})$ the associated analytic extension.
We recall that $\hat \A$ is the subspace of the algebraic dual of
$\A,$ consisting  of all functionals $ a \varphi $,
where  $a \in \A$. Here, $(a \f)(b)=\f(ba),$ and similarly,
 $(\varphi a)(b) = \f(ab), \ a,b \in \A.$ Since $\f a=\rho(a)\f$, we have
$\hat \A={\{
\varphi a\ | \ a \in \A \}}.$

A right-invariant positive linear functional $\hat{\psi}$ is
defined on $\hat{\A}$ by setting
$\hat{\psi}(\hat{a})=\varepsilon (a)$, for all $a\in \A$. Here
$\hat{a}=a\varphi.$ Since the
linear map, $\A\rightarrow\hat{\A}$, $a\mapsto\hat{a}$, is a
bijection (by faithfulness of  $\varphi$), the functional
$\hat{\psi}$ is well defined. Further, we have
$\hat{\psi} (\hat{b}^*\hat{a})=\varphi (b^*
a),$ for all $a,b\in \A.$

As shown in \cite{KuVD}, one may assume that the regular representation
$\hat \pi_{r}$ of $\hat \A$ also acts on $\H$.
Accordingly, we identify $\hat \A_{r}$ with
the C*-algebra generated by
 $\hat \pi_{r} (\hat \A)$ and set
$\hat \M = \hat \A_{r}''.$
A useful fact is that both $\M$ and $\hat \M$ act standardly on $\H.$

\vspace{1ex}
 We will  work  quite often with the ``opposite'' dual quantum group
 $(\hat\A_r,\hat\Delta_{r, \rop}) .$
Note that when we
add $\rop$ as a subscript to a co-product map, we  mean by this the
opposite co-product; that is, the one obtained after ``flipping'' the
original map. One reason for working with $(\hat\A_r,\hat\Delta_{r,
\rop})$ is that it corresponds to the dual of $(\A_r,\Delta_{r})$ as
defined in
\cite{KV}. Further, the multiplicative unitary associated to
$(\hat\A_r,\hat\Delta_{r,
\rop})$ is simply given by $\hat W = \Sigma W^* \Sigma,$
which fits
 with the usual notation for multiplicative unitaries and their
duals (cf.\ \cite{BS}).

\vspace{2ex}

We denote by $(\A_u,\Delta_u)$  the universal (locally compact) C*-algebraic
quantum
group associated to $(\A,\Delta)$, as introduced by J. Kustermans
\cite{Kus}. We recall
here some details of his construction.

The
 C$^{*}$-algebra $\A_u$  is the completion of $\A$ with respect to the
 C$^{*}$-norm $\| \cdot \|_{u} $ on $\A$ defined by
 \[ \| a \|_{u} = \sup \{ \| \Phi(a) \| \ |  \Phi \  \textnormal{is a
 $*$-homomorphism from $\A$ into some C$^{*}$-algebra}  \} \]
 (The non-trivial fact that this expression gives a well-defined norm on
$\A$
is shown in \cite{Kus}).
 The C*-algebra
  $\A_u$ has then the universal property that one may extend from
$\A$ to $\A_u$
  any $*$-homomorphism
 from $\A$ into
 some C$^{*}$-algebra.

 \vspace{2ex}

 The definition of $\Delta_{u}$ relies on the following
 proposition \cite[Proposition 3.8]{Kus}, which we restate here as
 we will need it in the sequel.

 \begin{proposition} \label{Kust}
     Consider   C*-algebras $C_{1}, C_{2}$ and $
 ^{*}$-homomorphisms $\phi_{1}$ from $\A$ into $M(C_{1})$ and
 $\phi_{2}$ from $\A$ into $M(C_{2})$ such that $\phi_{1}(\A) C_{1}$ is
 dense in $C_{1}$ and  $\phi_{2}(\A) C_{2}$ is
 dense in $C_{2}.$
 Then there exists a unique  $^{*}$-homomorphism $\phi$ from $\A$ into
$M(C_{1}
 \otimes C_{2})$ such that
 $$ ( \phi_{1}(a_{1}) \otimes \phi_{2}(a_{2})) \phi (a) =
  (\phi_{1} \odot \phi_{2}) ( (a_{1} \otimes a_{2}) \Delta(a) ) $$
  and
  $$ \phi (a)( \phi_{1}(a_{1}) \otimes \phi_{2}(a_{2}))  =
  (\phi_{1} \odot \phi_{2}) (\Delta(a) (a_{1} \otimes a_{2})) $$
  for every $a_{1}, a_{2} \in \A.$ We have moreover that $\phi(\A)
  (C_{1} \otimes C_{2})$ is
 dense in $C_{1} \otimes C_{2}.$

     \end{proposition}

     Now, let $\pi_{u}$ denote the identity
 mapping from $\A$ into $\A_u$. Hence, $\pi_{u}$ is an
 injective $^{*}$-homomorphism
 from $\A$ to $\A_u$ such that $\pi_{u}(\A) $ is dense in $\A_{u}$,
 so  $\pi_{u}(\A) \A_{u} $ is dense in $\A_{u}$ (as $\A^{2} = \A$).
     By applying the above proposition with $\phi_{1}=\phi_{2}= \pi_{u}$
     and exploiting the universal property of $\A_{u}$, one obtains
     that there exists a unique non-degenerate $^{*}$-homomorphism
     $\Delta_u : \A_u \to M(\A_u \otimes \A_u)$
  such  that
 \[ (\pi_{u} \odot \pi_{u})(x) \Delta_u (\pi_{u}(a)) =
 (\pi_{u} \odot \pi_{u})(x \Delta(a))   \]
 and
 \[ \Delta_u  (\pi_{u}(a)) (\pi_{u} \odot \pi_{u})(x) =
 (\pi_{u} \odot \pi_{u})(\Delta(a) x)   \]
 for all $a \in \A$ and $x \in \A \odot \A.$

  Being a $*$-homomorphism from $\A$ onto $\Complessi$, the co-unit
 $\varepsilon$ of $(\A, \Delta)$ extends to a
 $*$-homomorphism $\varepsilon_{u}$ from $\A_{u}$ onto
$\Complessi$, which is easily seen to satisfy the co-unit property
for $(\A_u,\Delta_u).$ Of course, we identify implicitly here $\A$
with its canonical copy $\pi_{u}(\A)$ inside $\A_{u}.$ Note that sometimes
we  add $u$ as an index
to denote the extension
to $\A_{u}$ of
a $^{*}$-homomorphism from $\A$ into some C$^*$-algebra, and sometimes
just use the same symbol to denote the extension when there is no
danger of confusion. For example, we get a canonical map $\pi_{r}:
\A_{u} \to \A_{r}$ which is the extension of $\pi_{r}:
\A \to \A_{r}.$

\vspace{1ex}
Let  now $(\A, \Delta)$ be an algebraic quantum group of
\emph{compact} type, that is, $\A$ has a
unit $I.$ It is immediate that $(\A_{r},\Delta_{r})$ is a
compact quantum group in the sense of Woronowicz \cite{Wo1,Wo2},
with Haar state $\varphi_{r}$  given by the restriction of the vector
state
$\omega_{\Lambda(I)}$ to
$\A_{r}.$
The  unique dense Hopf $*$-subalgebra \cite{BMT1} of $(\A_{r},\Delta_{r})$
may
 be identified with
 $(\A, \Delta, \varepsilon, S)$ ( via the Hopf $*$-algebra isomorphism
$\pi_{r}$).
 Using this identification, we may introduce the remarkable family
 $(f_{z})_{z \in \Complessi}$ of multiplicative linear functionals on $\A$
 constructed by Woronowicz (see \cite{Wo1, Wo2}).

 \vspace{1ex}
 Some of the properties of this family are:  $f_{0}= \e;$ $f_{z} *
 f_{z'} = f_{z+z'},$ where $\omega * \eta = (\omega \otimes \eta)
 \Delta, \ \
 \omega , \eta \in A';$ the maps $a \to f_{z} * a = (\iota \otimes f_{z})
 \Delta (a) $ and $ a \to a * f_{z} = (f_{z} \otimes \iota) \Delta (a) $ are
automorphisms
 of $\A;$ we have $f_{z}^* = f_{-\bar{z}}$ and $ f_{z} \circ S =
 f_{-z};$  we have $\f (a b) = \f ( b
 (f_{1}*a*f_{1}))$ and $S^2 (a) = f_{-1}* a * f_{1}.$

 \vspace{1ex} It follows from \cite[Theorem 2.12]{Ku} that
 $M(\hat \A),$ the multiplier algebra of $\A,$ may be concretely realized as
the
 subspace of the algebraic dual of $\A$ consisting of elements
 $\theta$ such that $(\theta \odot \iota) \Delta(a) $
 and $( \iota \odot \theta ) \Delta(a) $ belong to $\A$ for every
 $a \in \A.$ Hence, in the compact case, we have $f_{z} \in M(\hat \A) $ for
all $z.$

  \vspace{1ex} We also  recall that the following three conditions are
equivalent:

 \vspace{1ex}
 \ \ \ \ \ \ $\f$ is tracial; $f_{z} = \e $ for all $z \in \Complessi;$
$f_{1}= \e.$

 \vspace{1ex}
 The following description of algebraic quantum groups of
 \emph{discrete} type, that is, those which are dual to algebraic quantum
groups of
 compact type, will be useful. 

 \begin{proposition} \label{discrete} Let  $(\A, \Delta)$ be an algebraic quantum group of
compact type and
 let $(U^{\alpha})_{\alpha \in A}$ denote a complete set of
 pairwise inequivalent irreducible  unitary
 corepresentations of the compact quantum group $(\A_{r},\Delta_{r}).$
 Note that we have  $U^{\alpha} \in \A \odot M_{d_{\alpha}}(\Complessi)$
for some $d_{\alpha} < \infty$, when identifying $\A$ as the dense
Hopf $*$-algebra of $\A_{r}.$  Write each $U^{\alpha}$ as a matrix
$(u^{\alpha}_{ij})$ over $\A$
and recall that the set Span $ \{ \ (u^{\alpha}_{ij}) | 1 \leq i, j \leq
d_{\alpha}, \alpha \in A \ \}$ is a linear basis for $\A.$ Let $$M_{\alpha}=
\sum_{i=1}^{d_{\alpha}} f_{-1}(u^{\alpha}_{ii})
= \sum_{i=1}^{d_{\alpha}} f_{1}(u^{\alpha}_{ii})$$ denote the quantum
dimension of $U^{\alpha}.$

\vspace{1ex} \noindent
Further, set $\hat{\A}_{\alpha} = $ Span $ \{ \
(\hat{u}^{\alpha}_{ij}) | 1 \leq i, j \leq
d_{\alpha}\ \}$ and define
$$p_{\alpha}=  M_{\alpha} \sum_{i,j=1}^{d_{\alpha}} f_{1}(u^{\alpha}_{ji})
\hat{u}^{\alpha}_{ij}
\in \hat{\A}_{\alpha}.$$ Then

\vspace{1ex}
\noindent (1)
$$(i) \ \ \hat{u}^{\alpha}_{ij} \hat{u}^{\beta}_{kl} = \frac{\delta_{\alpha
\beta}}{M_{\alpha}} f_{-1}(u^{\alpha}_{kj} ) \hat{u}^{\alpha}_{il},$$
$$(ii) \ \ (\hat{u}^{\alpha}_{ij})^{*} = \hat{u}^{\alpha}_{ji},$$
$$(iii) \ \ \ p_{\alpha} \hat{u}^{\beta}_{kl} = \delta_{\alpha \beta}
\hat{u}^{\beta}_{kl} = \hat{u}^{\beta}_{kl} p_{\alpha},$$
where $1 \leq i, j \leq
d_{\alpha}, 1 \leq k, l \leq
d_{\beta}, \alpha, \beta \in A.$

\vspace{1ex}
\noindent
(2) The set $ \{ \ (\hat{u}^{\alpha}_{ij}) | 1 \leq i, j \leq
d_{\alpha}, \alpha \in A \ \}$ is a linear basis for $\hat{\A}_{\alpha}.$

\vspace{1ex}
\noindent
(3) Each $\hat{\A}_{\alpha}$ is a $*$-subalgebra
of $\hat \A,$ which is unital with unit $p_{\alpha}.$
As a $*$-algebra, $\hat{\A}_{\alpha}$ is isomorphic to the matrix
algebra
$M_{d_{\alpha}}(\Complessi).$

\vspace{1ex}
\noindent
(4) $\hat \A = \bigoplus_{\alpha} \hat{\A}_{\alpha}$ (algebraic direct
sum).

\vspace{1ex}
\noindent
(5) For each $\alpha \in A,$ let $Tr_{\alpha}$ denote the canonical trace
on $\hat{\A}_{\alpha} = M_{d_{\alpha}}(\Complessi)$
satisfying  $Tr_{\alpha}(p_{\alpha}) = d_{\alpha}.$ Then
$$ \hat \psi (x) = \oplus_{\alpha} Tr_{\alpha}(  p_{\alpha} x
f_{-1}), \ \ x \in \hat \A.$$

\end{proposition}

\begin{proof}  This result is  essentially known (see \cite[p. 
 722]{ER}), but we will need the explicit description presented here 
 in the sequel. We give a
proof for the sake of completeness. It relies  on the so-called
orthogonality relations for
the $U_{\alpha}$'s established in \cite{Wo1,Wo2}.

\vspace{1ex}
\noindent
(1) With obvious index notation, we have
$$ (\hat{u}^{\alpha}_{ij} \hat{u}^{\beta}_{kl})( (u^{\gamma}_{pq})^{*}) =
 \sum_{r} \f((u^{\gamma}_{pr})^{*}u^{\alpha}_{ij})
\f((u^{\gamma}_{rq})^{*}u^{\beta}_{kl})$$
$$ =(1/M_{\alpha})(1/M_{\beta}) \sum_{r}\delta_{\alpha \gamma}
f_{-1}(u^{\alpha}_{ip}) \delta_{rj} \delta_{\beta \gamma}
f_{-1}(u^{\beta}_{kr})
\delta_{ql}$$
$$=(1/M_{\alpha})(1/M_{\beta})\delta_{\alpha \gamma}\delta_{\beta
\gamma}\delta_{ql}
f_{-1}(u^{\alpha}_{ip})f_{-1}(u^{\beta}_{kj})$$
$$=(1/M^{2}_{\alpha})\delta_{\alpha \beta}\delta_{\alpha \gamma}\delta_{ql}
f_{-1}(u^{\alpha}_{ip})f_{-1}(u^{\alpha}_{kj})$$
$$=(1/M_{\alpha}) \delta_{\alpha \beta} f_{-1}(u^{\alpha}_{kj})
\f((u^{\gamma}_{pq})^{*}u^{\alpha}_{il})$$
$$= ((1/M_{\alpha}) \delta_{\alpha \beta} f_{-1}(u^{\alpha}_{kj})
\hat{u}^{\alpha}_{il}) ( (u^{\gamma}_{pq})^{*}),$$
and $(i)$ follows.
Concerning $(ii)$, we have
$$(\hat{u}^{\alpha}_{ij})^{*} =(S( u^{\alpha}_{ij})^*)^{\wedge} =
\hat{u}^{\alpha}_{ji},$$
using \cite[Lemma 7.14]{KuVD}. Further, using $(i)$, we get
$$ p_{\alpha} \hat{u}^{\beta}_{kl} = M_{\alpha} \sum_{i,j}
f_{1}(u^{\alpha}_{ji})
\hat{u}^{\alpha}_{ij} \hat{u}^{\beta}_{kl}$$
$$= M_{\alpha} (1/M_{\alpha}) \delta_{\alpha \beta} \sum_{i,j}
f_{1}(u^{\alpha}_{ji})
f_{-1}(u^{\alpha}_{kj}) \hat{u}^{\alpha}_{il} $$
$$=  \delta_{\alpha \beta} \sum_{i} \delta_{ik} \hat{u}^{\alpha}_{il}
= \delta_{\alpha \beta} \hat{u}^{\beta}_{kl}, $$
which is equal to $ \hat{u}^{\beta}_{kl} p_{\alpha}$ by a similar
computation. Hence, $(iii)$ is proved.

\vspace{1ex}
\noindent
(2) As $ \{ \ u^{\alpha}_{ij} | 1 \leq i, j \leq
d_{\alpha}, \alpha \in A \ \}$ is a linear basis for $ \A_{\alpha}$
and the map $a \to \hat a$ is a linear isomorphism between $\A$ and
$\hat{\A},$  (2) is clear.

\vspace{1ex}
\noindent
(3) The first sentence is an obvious consequence of (1). Now,
$\hat{\A}_{\alpha}$ may be seen as a finite dimensional C$^{*}$-algebra
(using the faithfulness of the $*$-homomorphism $\hat{\pi}_{r}$). Hence,
to show that $\hat{\A}_{\alpha}$ is isomorphic to
$M_{d_{\alpha}}(\Complessi),$ it is enough to show that
if $b = \sum_{i,j} b_{ij} \hat{u}^{\alpha}_{ij}, b_{ij} \in
\Complessi,$ is an element of
the center of $\hat{\A}^{\alpha},$ then $b= \lambda p_{\alpha},$
that is,
$$b_{ij}= \lambda f_{1}(\u_{ji}) \  \textnormal{for some} \ \lambda \in
\Complessi, 1 \leq i, j \leq d_{\alpha}.$$ Now, using (1)$(i)$, one
sees immediately that $\hat{u}^{\a}_{kl} b = b \hat{u}^{\a}_{kl}$ hold for
all $k$ and $ l$
if, and only if,
$$\sum_{i,j} b_{ij} f_{-1}(\u_{kj})\hat{u}^{\a}_{il} = \sum_{i,j} b_{ij}
f_{-1}(\u_{il})\hat{u}^{\a}_{kj}$$
for all $k, l,$ which in turn is equivalent to
$$(*) \ \ \ \ \delta_{rk} \sum_{i} b_{is} f_{-1}(\u_{il}) =
\delta_{sl}\sum_{j}b_{rj}f_{-1}(\u_{kj}), \ \ \forall k, l, r,
s.$$
We introduce now the two complex matrices $B=(b_{ij})$ and
$C= (c_{ij}),$ where $c_{ij}= f_{-1}(\u_{ji}).$ Then $(*)$ may be
rewritten as
$$(**) \ \ \ \ \delta_{rk} d_{ls} = \delta_{sl} e_{rk},  \ \ \forall \ k, l,
r,
s,$$
where $d_{ls} = \sum_{i} c_{li}b_{is}$ and $e_{rk}=
\sum_{j}b_{rj}c_{jk}.$ From $(**),$ we clearly get
$$ (BC)_{sl}= 0 = (CB)_{sl},  \  s \neq l \, ; \ \   (BC)_{ll} =
(BC)_{kk},$$
hence that $BC=CB$ is a  complex multiple $\lambda$ of the identity matrix.
But $C$ is invertible, with inverse $C^{-1}= (c'_{ij}),$ where
$c'_{ij}= f_{1}(\u_{ji}).$ Indeed, we have
$$ \sum_{j}c_{ij} c'_{jk} = \sum_{j} f_{-1}(\u_{ji})f_{1}(\u_{kj}) =
\e(\u_{ki}) = \delta_{ik}.$$
Therefore, we can conclude that $B = \lambda C^{-1},$ that is,
$b_{ij}= \lambda f_{1}(\u_{ji}).$ This etablishes (3). 

\vspace{1ex} \noindent (4) is an easy consequence of the previous 
assertions.

\vspace{1ex}
\noindent
(5) Fix now $\alpha \in A$ and define a linear functional $\tau$ on
$\hat{\A}_{\a}$ by $$\tau (x) = (1/M_{\a}) \hat\psi(xf_{1}), \ x \in
\hat{\A}_{\a}.$$
To show (4), a moment's thought makes it clear that
it is enough to show that $\tau = Tr_{\a}.$ Due to the uniqueness
property of $Tr_{\a},$ we only have to show that
$$(a) \ \ \tau \ \textnormal{is central} \,; \ \ \  (b)  \  \ \tau(p_{\a}) =
d_{\a}.$$
To show (a), we have to show $$ (a')  \ \ \ \hat\psi (xyf_{1}) =
\hat\psi(yxf_{1}) , \ \ x, y \in \hat{\A}^{\alpha}.$$
Now, let $\hat\rho $ denote the automorphism of $\hat{\A}$ satisfying
$\hat\psi(\hat a \hat b) = \hat\psi(\hat b \hat\rho(\hat a))$ for all
$a, b \in \A.$ Then we get $\hat\psi(yxf_{1}) = \hat\psi(xf_{1}
\hat\rho(y)),$
so $(a')$ follows if  $yf_{1}= f_{1}\hat\rho(y)$ hold for
all $y \in \hat{\A},$ that is, if $\hat \rho(y) = f_{-1} y f_{1}, \ y
\in \hat{\A}.$ This follows from Lemma \ref{roh}.

To show $(b),$ we first observe that
$$ (\hat{u}^{\a}_{ij} f_{1})(u^{\beta}_{kl})^*) = \sum_{r}
\f((u^{\beta}_{kr})^* \u_{ij}) f_{1}((u^{\beta}_{rl})^*)$$
$$= \delta_{\a \beta} (1/M_{\a}) \sum_{k}f_{1}((u^{\beta}_{rl})^*)
\delta_{rj}
f_{-1}(\u_{ik}) = \delta_{\a \beta} (1/M_{\a})
f_{1}((u^{\beta}_{jl})^*)f_{-1}(\u_{ik}).$$
Using this, we show that $p_{\a}f_{1}= M_{\a} \sum_{i}\hat{u}^{\a}_{ii}.$
Indeed, we have
$$(p_{\a}f_{1})((u^{\beta}_{kl})^*)  =
M_{\a} (\sum_{i,j}f_{1}(\u_{ji})\hat{\u}_{ij}f_{1})
((u^{\beta}_{kl})^*)$$
$$= \delta_{\a \beta}
\sum_{i,j}f_{1}(\u_{ji})f_{1}((u^{\beta}_{jl})^*)f_{-1}(\u_{ik})$$
$$= \delta_{\a \beta} \sum_{j} \delta_{jk}f_{1}((u^{\beta}_{jl})^*)
= \delta_{\a \beta} f_{1}((u^{\beta}_{kl})^*)$$
$$=\delta_{\a \beta} f_{1}(S(u^{\beta}_{lk})) = \delta_{\a \beta}
f_{-1}(u^{\beta}_{lk}),$$
while
$$ M_{\a} \sum_{i}\hat{u}^{\a}_{ii}((u^{\beta}_{kl})^*) = M_{\a}
\sum_{i} \f((u^{\beta}_{kl})^* \u_{ii}) $$
$$= \delta_{\a \beta} \sum_{i}\delta_{li} f_{-1}(u^{\beta}_{ik})
=\delta_{\a \beta} f_{-1}(u^{\beta}_{lk}).$$
But then we get
$$ \tau( p_{\a}) = (1/M_{\a}) \hat\psi(p_{\a}f_{1}) =
\sum_{i}\hat\psi(\hat{u}^{\a}_{ii}) $$
$$= \sum_{i} \e ( \u_{ii}) = d_{\a},$$
and $(b)$ is shown. This finishes the proof of (5) and of the
proposition.
\end{proof}

\begin{lemma} \label{roh}
 Let  $(\A, \Delta)$ be an algebraic quantum group of
compact type. Let $\hat\rho $ denote the automorphism of $\hat{\A}$
satisfying
$\hat\psi(\hat a \hat b) = \hat\psi(\hat b \hat\rho(\hat a))$ for all
$a, b \in \A.$ Then $$\hat\rho ( \hat a) = f_{-1} \hat a f_{1}, \ a
\in \A.$$
\end{lemma}
\begin{proof}
   Being of compact type, $(\A, \Delta)$ is unimodular, that is,
  the  modular element $\delta$ of $\A$ is trivial. Hence, it follows from
   \cite[Lemma 2.2]{BMT3} that $\hat\rho ( \hat a) =
   (S^2(a))^{\wedge}$ for all $  a \in \A.$ Therefore, we have
   $$\hat\rho ( \hat a) (b) = (S^2(a))^{\wedge}(b) = \f (bS^2(a))$$
   $$= \f S^{-2} ( b S^2(a)) = \f ( S^{-2}(b)a) = \hat a(S^{-2}(b))$$
   $$= \hat a(f_{1}*b*f_{-1}) = \hat a ( (f_{-1} \odot \iota \odot
   f_{1})(\Delta \odot \iota) \Delta(b))$$
   $$= (f_{-1}\odot \hat a \odot f_{1}) (\Delta \odot \iota)
   \Delta(b) = (f_{-1} \hat a f_{1}) (b)$$
   for all $a, b \in \A,$ which proves the assertion.
   \end{proof}

\section{Categorical interlude}
In this section we introduce
Rep$(\A_u,\Delta_u),$  Rep$(\A,\Delta)$ and Corep$(\hat
\A_r,\hat\Delta_{r,\rop})$
as concrete tensor
$C^*$-categories, and describe explicitly  isomorphisms
(of tensor $C^*$-categories) between them.
We also briefly mention some related categories. Finally,  we
establish the absorbtion property for the regular representation with
respect to tensor product.

\medskip
\subsection{$C^*$-tensor categories associated with algebraic quantum
groups}
We refer to \cite{GLR,LR} for terminology concerning tensor $C^*$-tensor
categories.
Let $(\A,\Delta)$ denote an algebraic quantum group. We begin with the
category $\R$ = Rep$(\A_u,\Delta_u)$
and explain how it may be organized as a concrete tensor
$C^*$-category with irreducible unit.
The objects in $\R$ are the $*$-representations $\pi$ of $\A_u$
acting on a Hilbert space $H_\pi$
satisfying the non-degenerateness (denseness) condition
$\overline{\pi(\A_u)\H_\pi} = \H_\pi.$
The family of arrows
(or morphisms)
between two objects $\pi$ and $\pi'$ is given by
$$Mor(\pi,\pi') =
\{T \in B(\H_\pi,\H_{\pi'}) \ | \  T\pi(a)=\pi'(a)T,
\ \forall a \in \A\} \ .$$
The element $1_{\pi} \in$
$Mor(\pi,\pi)$ is given by the identity on $\H_{\pi}$.
The adjoint of an element $T \in Mor(\pi,\pi')$ is given by its Hilbert
space adjoint $T^* \in$ $Mor(\pi',\pi),$ so  we  clearly have $\| T^{*}\ T
\| =
\|T\| ^2$.
The tensor product $\pi \times \pi'$ of two objects $\pi$ and $\pi'$ is
defined as
$\pi \times \pi'=(\pi \otimes \pi')\Delta_{u},$
while
on the arrows we have
the usual tensor product of operators.
The unit in the tensor category is given by $\varepsilon_{u}$. Note that
this unit is irreducible, since
$Mor(\varepsilon_{u},\varepsilon_{u})= \Complessi$. It is clear  that
this category has natural subobjects and that one may form direct sums
in an obvious way.

\medskip
Next, we introduce the closely related category $\R_{alg}$ =
Rep$(\A,\Delta).$
The objects in $\R_{alg}$ are now the $*$-representations $\pi$ of $\A$
acting on a Hilbert space $H_\pi$
satisfying the non-degenerateness (denseness) condition
$\overline{\pi(\A)\H_\pi} = \H_\pi.$ Arrows and adjoints are defined
in a similar way as above.  To define the tensor product of objects, we have
to
appeal to Proposition \ref{Kust}. Let  $\phi_{1}$ and $\phi_{2}$ be
objects in $\R_{alg}$, and consider  $\phi_{1}$ (resp. $\phi_{2}$) as
a $^{*}$-homomorphism from $\A$ into $M(B_{0}(\H_{\phi_{1}}))$
(resp. $M(B_{0}(\H_{\phi_{2}})).$ As $\phi_{1}$ and $\phi_{2}$ are
non-degenerate (by assumption), the proposition applies
and produces a unique
$^{*}$-homomorphism $ \phi_{1} \times \phi_{2}$ from $\A$ into $
M(B_{0}(\H_{\phi_{1}})\otimes B_{0}(\H_{\phi_{2}}) ) = M (
B_{0}(\H_{\phi_{1}} \otimes  \H_{\phi_{2}})) = B (\H_{\phi_{1}} \otimes
\H_{\phi_{2}}) $
such that
 $$ ( \phi_{1}(a_{1}) \otimes \phi_{2}(a_{2})) (\phi_{1} \times \phi_{2})
(a) =
  (\phi_{1} \odot \phi_{2}) ( (a_{1} \otimes a_{2}) \Delta(a) ) $$
  and
  $$ (\phi_{1} \times \phi_{2})  (a)( \phi_{1}(a_{1}) \otimes
\phi_{2}(a_{2}))  =
  (\phi_{1} \odot \phi_{2}) (\Delta(a) (a_{1} \otimes a_{2})) $$
  for every $a_{1}, a_{2} \in \A.$ Moreover, we have that $(\phi_{1} \times
  \phi_{2}) (\A)
  (B_{0}(\H_{\phi_{1}})\otimes B_{0}(\H_{\phi_{2}}))$ is
 dense in $B_{0}(\H_{\phi_{1}})\otimes B_{0}(\H_{\phi_{2}}).$ It
 follows easily that   $\phi_{1} \times \phi_{2}$, when regarded as a
 $^{*}$-homomorphism from $\A$ into $B(\H_{\phi_{1}} \otimes
 \H_{\phi_{2}})$, satisfies the non-degeneracy (denseness) condition
 required for qualifying it as an object in $\R_{alg}.$
 Finally, the unit in the tensor category is of course $\varepsilon.$

\smallskip
Not surprisingly, the following proposition holds:

\begin{proposition} Define $P: \R \to \R_{alg}$ on objects by
    $P(\pi) = \pi \circ \pi_u, $ and let $P$ act trivially on arrows. Then
$P$ is an isomorphism of tensor $C^*$-categories.
\end{proposition}

\begin{proof} The only non-trivial fact in this assertion is perhaps to show
    that $P$ preserves tensor products.

    Let $\pi_{1}, \pi_{2}$ be objects in $\R,$ and set $\phi_{1}
    =P(\pi_{1}), \phi_{2} = P( \pi_{2}), \phi = P(\pi_{1} \times
    \pi_{2}).$ We have to show that $ \phi = \phi_{1} \times
    \phi_{2}.$

    Now, let $a, a_{1}, a_{2} \in \A.$ Then we have

    $$ \phi (a) ( \phi_{1} (a_{1}) \otimes \phi_{2}(a_{2}) )$$
    $$ = (\pi_{1} \times \pi_{2} ) \pi_{u}(a) (  \pi_{1} \pi_{u}(a_{1})
    \otimes \pi_{2}\pi_{u}(a_{2}) )$$
    $$= (\pi_{1} \otimes \pi_{2} ) (\Delta_{u}\pi_{u}(a))  (\pi_{1}
    \otimes \pi_{2}) (\pi_{u}(a_{1}) \otimes \pi_{u}(a_{2})) $$
    $$= (\pi_{1} \otimes  \pi_{2} ) ( \Delta_{u}\pi_{u}(a)
    ( \pi_{u}(a_{1}) \otimes \pi_{u}(a_{2})) )$$
    $$= (\pi_{1} \odot \pi_{2} )  (\pi_{u} \odot \pi_{u}) ( \Delta
    (a) ( a_{1} \otimes a_{2}) )$$
    $$ = ( \phi_{1} \odot \phi_{2} ) ( \Delta (a) ( a_{1} \otimes
    a_{2})).$$
    In the same way, one shows that
    $$  ( \phi_{1} (a_{1}) \otimes \phi_{2}(a_{2}) ) \phi (a)
    = ( \phi_{1} \odot \phi_{2} ) ( ( a_{1} \otimes
    a_{2}) \Delta (a)).$$
    From the uniqueness property of $\phi_{1} \times \phi_{2},$, we
    may then conclude that $\phi = \phi_{1} \times \phi_{2},$ as
    desired.

\end{proof}

\medskip

Another category $\C=$ Corep$(\hat\A_r,\hat\Delta_{r,\rop})$ which may
be organized as a concrete tensor
$C^*$-category is defined as follows. The objects in $\C$ consist of
unitary elements
$U$ lying in $(M(\hat\A_r \otimes B_0(\H_{U})))$ for some Hilbert
space $\H_{U}$ and  satisfying
the corepresentation property
$$(\hat\Delta_{r, \rop} \otimes \iota)U = U_{13} \, U_{23} .$$

\noindent For objects $U$ and $V$  in $\C,$ we set
$$Mor(U,V) =
\{T \in B(\H_U,\H_V) \ | \
T (\omega \bar\otimes \iota) (U )= (\omega \bar\otimes \iota) (V) T \
, \ \forall \omega \in \hat{\M}_* \}. $$ The element $1_{U} \in$
$Mor(U,U)$ is given by the identity on $\H_{U}.$
The adjoint of an element $T \in Mor
(U,V)$ is given by its Hilbert
space adjoint $T^* \in Mor(V,U)$, so  we  have $\| T^{*}\ T \| =
\|T\| ^2.$
The monoidal structure on the objects is determined by setting
$$U \times V = V_{13} \, U_{12} \in M(\hat \A_{r} \otimes B_0(\H_U)
\otimes B_0(\H_V)) =M(\hat \A_{r} \otimes B_0(\H_U \otimes
\H_{V})).$$
(Clearly $U \times V$ is unitary. Moreover,
\[(\hat\Delta_{r, \rop} \otimes \iota) (U \times V) =
(\hat\Delta_{r, \rop} \otimes \iota \otimes \iota )( V_{13} \, U_{12} ) \]
\[ = (\hat\Delta_{r, \rop} \otimes \iota \otimes \iota )( V_{13})
(\hat\Delta_{r, \rop} \otimes \iota \otimes \iota )( U_{12}) \]
\[ = V_{14}  \,  V_{24}  \,  U_{13}  \, U_{23} = V_{14}  \,  U_{13}  \,
V_{24}
 \, U_{23} = ( U \times V )_{13} ( U \times V )_{23}, \]
 so $U \times V \in \C.$)
The reason for ``reversing'' the
 ``natural'' tensor product
will be evident from Theorem \ref{cateq}.
Again, on the arrows, we just set $T \times S = T \otimes S$.

The unit in the tensor category is given by $I \otimes 1 \in M(\hat\A_r
\otimes
\Complessi) $, where $I$
denotes the unit of $M(\hat\A_{r})$.
As we clearly have $Mor(I \otimes 1,I \otimes 1)= \Complessi,$
this unit is irreducible.

\vspace{1ex}
For objects $U,V$ in $ \C,$ we say that $U$ is (unitarily) \emph{equivalent}
to $V$,
and write $U \simeq V$,
whenever there exists a unitary $T \in Mor(U,V).$

\vspace{2ex}
One may clearly introduce several  related tensor
 C$^{*}$-categories, such as \\
Rep$(\A_{u},\Delta_{u,\rop})$,
Corep$(\A_{u},\Delta_{u}),$
Corep$(\A_r,\Delta_{r})$ and
Corep$(\A_r,\Delta_{r,\rop})$, along the same lines. Note that in the
sequel we  always refer
to the monoidal structure defined in the same way as above.
Nevertheless, it should be noted
that there are relations
between the possible choices of monoidal stucture in these tensor
categories. For example, if one considers Corep$(\A_{r},\Delta_{r})$ with
monoidal structure given by $X \bar{\times} Y = X_{12}Y_{13}$ and
Corep$(\A_{r},\Delta_{r,\rop})$ with $U \times V = V_{13} U_{12},$
then it
is  easy
to check that the map $X \to X^*$ gives an
isomorphism between these two tensor categories (acting trivially on
arrows).

\subsection{From corepresentations to representations and back}
\vspace{2ex}
We now recall some results from \cite{BMT3}.
First, there exists an injective,
not necessarily $*$-preserving,
homomorphism
$Q_r : \A \to \hat\A_r^* $
determined by
$$Q_r(a)(\hat\pi_r(\hat{b}))=\hat{b}(S^{-1}(a))=\varphi(S^{-1}(a)b),
\quad \forall a,b \in \A.$$
In fact, there  exists an injective homomorphism
$Q : \A
\to
\hat \M_{*}$ satisfying \\ $Q(a)_{|\hat \A_{r}}= Q_{r}(a)$ for all $a \in
\A,$ and
such that $Q(\A)$ is norm-dense in $\hat \M_{*}$.
For all $a \in \A$, we have  $Q(a) =
\omega_{\Lambda(a),\Lambda(c)}$ (restricted to $\hat \M$),
where $c \in \A$ is
chosen  such that $\hat c \ \widehat{S(a^{*})} =
\widehat{S(a^{*})}$ (such a choice is always possible).

We will need
the following lemma.

\begin{lemma}Let $a \in \A $ and assume that $c \in \A$ satisfies
    $\hat c \ \widehat{S(a^{*})} = \widehat{S(a^{*})}.$

    Then we have

    $$\Delta(a) = ( \iota \odot \iota \odot \varphi) ((I
    \otimes I \otimes c^{*} ) ( \iota \odot \Delta) \Delta (a)  ).$$

\end{lemma}

\begin{proof} Recall from \cite{BMT3} that

    $$\hat c \ \widehat{S(a^{*})} = \sum_{i} \varphi(q_{i}) \hat
    p_{i},$$

    \hspace{-2em}where $ c \otimes S(a^{*}) = \sum_{i} \Delta(p_{i}) (q_{i}
    \otimes I)$ for some $p_{i}, q_{i} \in \A, (i=1\ldots n).$

     \vspace{2ex}
    Note
    also that the inversion formula from \cite[p. 1024]{KuVD} gives
    then

    $$\sum_{i} q_{i} \otimes p_{i} = ((S^{-1} \odot \iota) \Delta
    S(a^{*})) ( c \otimes I). $$

    Hence, we have

    $$ (\hat c \ \widehat{S(a^{*})})(b) = \sum_{i} \varphi(q_{i}) \hat
    p_{i} (b)$$

    $$= \sum_{i} \varphi(q_{i}) \varphi(b p_{i}) = (\varphi \odot
    \varphi) ( (I \otimes b) ( \sum_{i} q_{i} \otimes p_{i}))$$

    $$ = (\varphi b) ( (\varphi \odot  \iota) ( \sum_{i} q_{i} \otimes
    p_{i})  ) = (\varphi b) (  \, (\varphi \odot \iota) ( ((S^{-1} \odot
\iota) \Delta
    S(a^{*})) ( c \otimes I) ) \, ) $$

    while

    $$ \widehat{S(a^{*})}(b) = \varphi(b S(a^{*}) = (\varphi b)
    (S(a^{*}))$$
    for all $b \in \A.$

    From the assumption and the fact that
    $\{\varphi  b \ | \ b \in A \}$ separates points in $\A $
    (since $\varphi $ is faithful on $\A$), we get

    $$ S(a^{*}) = (\varphi \odot \iota) ( ((S^{-1} \odot \iota) \Delta
    S(a^{*})) ( c \otimes I) ) = (c \varphi \otimes \iota ) (S^{-1} \odot
\iota) \Delta
    S(a^{*})$$
    $$ = (\iota \odot c \varphi) ( S \odot \iota) \Delta(a^{*})
    = S (  \iota \odot c \varphi) \Delta (a^{*}),$$

   where we have used that $\Delta S = \chi (S \odot S) \Delta.$
   Therefore, we have

    $$a^{*} = ( \iota \odot c \varphi) \Delta (a^{*}) =
    ( \iota \odot  \varphi)( \Delta (a^{*}) (I \otimes c)),$$

    so $$ a = ( \iota \odot  \varphi)( (I \otimes c^{*}) \Delta (a) ),$$
    as $\iota \odot \varphi$ is $*$-preserving.

     \vspace{2ex}
    This implies that

    $$\Delta(a) = \Delta ( \iota \odot  \varphi)( (I \otimes c^{*}) \Delta
    (a) )$$

    $$ = (\iota \odot \iota \odot \varphi) (\Delta \odot \iota) ( (I
    \otimes c^{*}) \Delta(a) ) $$

    $$= (\iota \odot \iota \odot \varphi)
    ((I \otimes I \otimes c^{*})( \Delta \odot \iota) \Delta(a))$$

    $$= (\iota \odot \iota \odot \varphi) ((I \otimes I \otimes
    c^{*}) ( \iota \odot \Delta) \Delta(a) ) ,$$
    as asserted.

\end{proof}

The following result teaches us that $\C$ and $\R$ are isomorphic as
tensor C*-categories.

\begin{theorem}
\label{cateq}
Define $F: \C \to \R$  on objects by $F(U)=\pi_U$,
where $\pi_U$ is determined by
$$\pi_U(a) = (Q_r(a) \otimes \iota)U, \ \forall a \in \A,$$
and let $F$ act identically on arrows.

Define $G: \R \to \C$ on objects by $G(\pi)=U_\pi$,
where $U_\pi$ is determined by
$$U_\pi(\Lambda(a) \otimes \pi(b)v) = \sum_{i=1}^{n} \Lambda(a_i) \otimes
\pi(b_i) v$$ for all $ a,b \in \A, v \in H_{\pi},$ and the $a_{i}$'s
and $b_{i}$'s are elements in $\A$ chosen as to satisfy
$\Delta(a) (b \otimes I) = \sum_{i=1}^{n} b_{i} \otimes a_{i}. $
Let  $G$ act identically on arrows.

 \vspace{2ex}
Then $F$ and $G$ are  covariant monoidal (tensor preserving) functors
which are adjoint- and unit-preserving, and satisfy $GF={\rm id}$, $FG =
{\rm id}$.

 \vspace{2ex}
We also have  $F(\hat{W}) = \pi_r$ and $F(I \otimes 1) =
\varepsilon_{u}$.
\end{theorem}

\begin{proof} The fact that $F$ and $G$ are well defined on objects is
established
    in \cite{BMT3}, where it is also shown that
    $GF={\rm id}$, $FG = {\rm id}$, $F(\hat{W}) = \pi_r$ and $F(I \otimes 1)
= \varepsilon_{u}$.

     \vspace{2ex}
    We now check that  $F$ and $G$ are well defined on
    arrows.
    Let $U, V $ be unitary corepresentations of
    $(\hat\A_r,\hat\Delta_{r,\rop}),$ and let $T \in Mor(U,V)$. Then,
    for all $a \in \A$,  we have
    $$ T \ ( Q_{r}(a) \otimes \iota ) (U) =  ( Q_{r}(a) \otimes \iota )
    (V) \ T, $$
    hence, $T \ \pi_{U}(a) = \pi_{V}(a) \ T$. As $\A$ is norm-dense
    in $\A_{u}$, we get $T \pi_{U}(x)= \pi_{V}(x) T$ for all $x \in
    \A_{u}$, that is $T \in Mor(\pi_{U},\pi_{V})$. Conversely, if
    $T \in Mor(\pi_{U},\pi_{V})$, then, using that $Q(\A) $ is
    norm-dense in $\hat {\M} _{*}$, one readily sees that $ T \in
Mor(U,V).$

    \vspace{2ex}
    It is obvious that  $F$ and $G$ are adjoint- and
    unit-preserving. To show that $F$ and $G$
    are  monoidal, that is, preserve tensor products, it is
    enough to show that $F(U \times V) = F(U) \times F(V) $, that is,
    $\pi_{U \times V} = \pi_{U} \times \pi_{V}$, where $U$ and $V$ are
    unitary corepresentations of $(\hat\A_r,\hat\Delta_{r,\rop}).$

    \vspace{2ex}

    Let $a,b,f \in A, \xi \in \H_{U}, \eta \in \H_{V}. $ Then we have
    $$[(\pi_{U} \times \pi_{V}) (a)] ( \pi_{U}(b) \xi \otimes \pi_{V}(f)
\eta )
    = [(\pi_{U} \times \pi_{V}) (a)] ( \pi_{U}(b) \ \odot
    \pi_{V}(f)  ) (\xi \otimes \eta)$$
    \hspace{3em} $= (\pi_{U} \odot \pi_{V}) ( \Delta(a) (b \otimes f) ) (\xi
\otimes
    \eta).$

    \vspace{2ex}

    On the other hand, choose $c \in A$ such that  $\hat c \
\widehat{S(a^{*})}
    = \widehat{S(a^{*})}$, so we have $Q_{r}(a)=
    \omega_{\Lambda(a),\Lambda(c)}$ (restricted to $\hat \A_{r}$). Then
    $$\pi_{U \times V} (a) = (Q_{r}(a) \otimes \iota \otimes \iota) (
    U \times V) = (\omega_{\Lambda(a),\Lambda(c)} \otimes \iota \otimes
\iota)
    (V_{13} U_{12}).$$
    Now, choose $h_{k}, g_{k} \in \A, (k=1\ldots n)$, such that
    $\Delta(a) (b\otimes I) = \sum_{k} h_{k} \otimes g_{k}.$

    \vspace{2ex}
    Further,
    for each $k,$ choose $f_{l} ,  g_{k}^{l} \in \A \  (l=1\ldots
    n_{k}),$ such that

    $$\Delta(g_{k}) (f\otimes I) = \sum_{l} f_{l} \otimes
    g_{k}^l.$$

    We  then get

     \vspace{2ex}

    $ ( \ [(\pi_{U \times V}) (a)] ( \pi_{U}(b) \xi \otimes \pi_{V}(f)
    \eta ) \ , \ \xi' \otimes \eta' )$

    \vspace{2ex}
    $= (V_{13} U_{12} (\Lambda(a) \otimes \pi_{U}(b) \xi \otimes \pi_{V}(f)
    \eta ), \Lambda(c) \otimes \xi' \otimes \eta')$

    $$= \sum_{k} ( V_{13} (\Lambda(g_{k}) \otimes \pi_{U}(h_{k}) \xi \otimes
\pi_{V}(f)
    \eta , \Lambda(c) \otimes \xi' \otimes \eta') \ \
    \textnormal{(using the relation between} \ U  \ \textnormal{and}
    \ \pi_{U})$$
    $$= \sum_{k,l} (  \Lambda(g_{k}^l) \otimes \pi_{U}(h_{k}) \xi \otimes
\pi_{V}(f_{l})
    \eta ), \Lambda(c) \otimes \xi' \otimes \eta') \ \ \
    \textnormal{(using the relation between} \ V  \ \textnormal{and}
    \ \pi_{V})$$

     \hspace{-2em}
$= \sum_{k,l} \varphi(c^{*}g_{k}^l) \ ( (\pi_{U} \odot
\pi_{V})
    (h_{k} \otimes f_{l}) (\xi \otimes \eta) \ , \ \xi' \otimes
    \eta'),$

    \vspace{2ex}
    \hspace{-2em}for all $\xi' \in \H_{U}, \eta' \in \H_{V}$.

    \vspace{2ex}
    Therefore, it suffices to
    show that
    $$\Delta(a) (b \otimes f) = \sum_{k,l} \varphi(c^{*}g_{k}^l) (h_{k}
    \otimes f_{l}).$$
    Now, we have
    $$ \Delta(a) (b \otimes f) =\Delta(a) (b \otimes I)(I \otimes f),$$
    while
    $$ \sum_{k,l} \varphi(c^{*}g_{k}^l) (h_{k}
    \otimes f_{l}) = \sum_{k} (\iota \odot \iota \odot \varphi c^{*} ) (
h_{k}
    \otimes ( \sum_{l} f_{l} \otimes g_{k}^l ) )$$
    $$   = \sum_{k} (\iota \odot \iota \odot \varphi c^{*} ) ( h_{k}
    \otimes \Delta(g_{k}) ( f \otimes I ) )  )$$
    $$=   \sum_{k} h_{k} \otimes  \ ( \iota \odot \varphi ) ( (I \otimes
    c^{*}) \Delta(g_{k}) (f \otimes I))  $$
   $$ =  (\sum_{k} h_{k} \otimes   ( \iota \odot \varphi ) ( (I \otimes
    c^{*}) \Delta(g_{k}) ) ) (I \otimes f).$$

    Hence, this reduces to showing

    $$ \Delta(a) (b \otimes I) = \sum_{k} h_{k} \otimes ( \iota \odot
\varphi ) ( (I \otimes
    c^{*}) \Delta(g_{k}) ).$$

    Now, using the previous lemma, we have

    $$\Delta(a) (b \otimes I) = (( \iota \odot \iota \odot \varphi) ((I
    \otimes I \otimes c^{*} ) ( \iota \odot \Delta) \Delta (a) ) )(b
    \otimes I)$$
    $$ =( \iota \odot \iota \odot \varphi) ((I
    \otimes I \otimes c^{*} ) ( \iota \odot \Delta) \Delta (a)  (b
    \otimes I \otimes I) \ )$$
    $$= ( \iota \odot \iota \odot \varphi) ((I
    \otimes I \otimes c^{*} ) ( \iota \odot \Delta) (\Delta (a)(b
    \otimes I )))$$
    $$= ( \iota \odot \iota \odot \varphi) ((I
    \otimes I \otimes c^{*} ) ( \iota \odot \Delta) ( \sum_{k} h_{k}
    \otimes g_{k})) $$
    $$= ( \iota \odot \iota \odot \varphi) ((I
    \otimes I \otimes c^{*} ) ( \sum_{k} h_{k} \otimes \Delta
    (g_{k})))$$
    $$=( \iota \odot \iota \odot \varphi) ( \sum_{k} h_{k} \otimes
    (I \otimes c^{*}) \Delta(g_{k}) \ )$$
$$= \sum_{k} h_{k} \otimes  ( \iota \odot \varphi )
    ( (I \otimes c^{*}) \Delta(g_{k}) ),$$
    which finishes the proof.

\end{proof}

We may dualize this result by using Pontryagin's duality for algebraic
quantum groups \cite{VD2,KuVD}.  Attached to $(\hat \A, \hat \Delta),$ we
can first associate an injective homomorphism $\hat Q_{r} : \hat \A \to
\hat{\hat {\A_{r}^*}} \simeq \A^*_{r}$ which is determined by
$$ \hat Q_{r} (\hat a) (\pi_{r}(b)) = \hat a (S(b)) = \varphi (S(b)a)$$
for all $a, b \in \A.$ Then we get a functor $\hat F : $
Corep$(\A_{r},\Delta_{r}) \to$ Rep$(\hat \A_{u}, \hat{\Delta}_{u,\rop})$
determined on objects by
$$ \hat F (U) ( \hat a ) = (\hat Q_{r} (\hat a) \otimes \iota) U ,
\ \hat a \in \hat \A,$$ and acting trivially on arrows,
which is an isomorphism of tensor C*-categories and satisfies $F(
W ) = \hat \pi_{r}, \hat F ( I \otimes 1) = \hat \varepsilon_{u}.$ We will
write $\hat \pi_{U}$ for $\hat F (U)$ in the sequel.

\subsection{The absorbing property for $\pi_r$ and $\hat W$}

We show  that $\pi_r$ and $\hat W$ have an absorbing property with
respect to tensoring, which is analogous to Fell's classical result
for the regular representation of a group \cite{Fell}.

\vspace{2ex}

\begin{proposition} \label{absor}
Let $U$ be a unitary corepresentation of $(\hat\A_r,\hat\Delta_{r,\rop})$
and $I_U = I \otimes I_{\H_U}$
be the trivial unitary corepresentation of $(\hat\A_r,\hat\Delta_{r,\rop}).$
Then
$U \times \hat{W}$ and $I_U \times \hat{W}$ are equivalent objects in $\C$.
\end{proposition}
\begin{proof}
Set
$T=\chi U^*
\in M(B_0(\H_U) \otimes \hat\A_r) \subset B(\H_U \otimes \H).$
It suffices to check that this unitary
satisfies the relation
$$T (\omega \bar\otimes \iota) (U \times \hat{W})
= (\omega \bar\otimes \iota) (I_U \times \hat{W}) T,
\quad \forall \omega \in \hat \M_* \ .$$
After some manipulations, this relation reduces to
$$(\omega \bar\otimes \iota \bar\otimes \iota) (U^*_{32} \hat{W}_{13}
U_{12})
= (\omega \bar\otimes \iota \bar\otimes \iota) (\hat{W}_{13} U^*_{32}).$$
Now, using the fact that $\hat W$ is a multiplicative unitary,
we get
$$ \hat{W}_{13}U_{12}
= U_{32}\hat{W}_{13}U^*_{32}
\in M(\hat\A_r \otimes B_0(\H_U \otimes \H)).$$
Thus, we have
$U^*_{32} \hat{W}_{13} U_{12} = \hat{W}_{13} U^*_{32}$,
and the result clearly follows.
\end{proof}

Combining this result with Theorem \ref{cateq},
one gets at once that
$\pi \times \pi_r$ is equivalent to
$I_{\pi} \times \pi_r$ for every $\pi \in $ Rep$(\A_u,\Delta_u)$,
where $I_\pi \in$ Rep$(\A_u,\Delta_u)$ is given by
$I_\pi(a) = \varepsilon_u(a) I_{\H_\pi}, \ \forall a \in \A_u$.

\vspace{2ex}
By duality, we also have a similar result for $\hat{\pi}_{r}$ and $W.$

\section{Conjugate and  Hilbert-Schmidt corepresentations}
In this section, we define the conjugate and the Hilbert-Schmidt
corepresentations associated with a unitary
corepresentation. Such objects play an important role in the classical
representation theory for groups and we will need these concepts in
later sections.

\subsection{Conjugate corepresentations}
Let $(\A,\Delta)$ be an algebraic quantum group
and $U$ be a unitary corepresentation of $(\A_r,\Delta_r).$
Let $\overline{\H}_U$ be any  Hilbert space  such that there exists an
anti-unitary map  $J : \H_U \to \overline{\H}_U.$
Define then $j: B(\H_U) \to B(\overline{\H}_U)$ by
$j(x)= Jx^*J^*, \ \forall x \in B(\H_U)$.
Then $j$ is linear, unital, normal, isometric,
$*$-preserving and anti-multiplicative,
with inverse $j^{-1}(\overline{x}) = J^* \overline{x}^* J, \
 \overline{x} \in B(\overline{\H}_U)$.
Note that $j(B_0(\H_U))
 = B_0(\overline{\H}_U).$
 
\vspace{1ex}
We may then define
$$\overline{U}
= (R \otimes j)U \in M(\A_r \otimes B_0(\overline{\H}_U)).$$

\begin{proposition}
 $\overline{U}$ is a unitary corepresentation of $(\A_r,\Delta_r),$
with $\H_{\overline{U}} = \overline{\H}_U.$
\end{proposition}

\begin{proof} We have
$$\overline{U}^* \overline{U}
= ((R \otimes j)U)^* ((R \otimes j)U)
= (R \otimes j)(UU^*)
= (R \otimes j)(I)= I_{M(\A_r)} \otimes I_{\H_U} \ ,$$
and similarly
$$\overline{U} \overline{U}^*
= I_{M(\A_r)} \otimes I_{\overline{\H}_U} \ .$$
Furthermore,
$$(\Delta_r \otimes \iota)\overline{U}
= (\Delta_r R \otimes j)U
= (\chi(R \otimes R) \otimes j)(\Delta_r \otimes \iota)U
= (\chi(R \otimes R) \otimes j)(U_{13}U_{23}) $$
$$= (\chi(R \otimes R) \otimes j)(U_{23})
(\chi(R \otimes R) \otimes j)(U_{13})
= ((R \otimes j)U)_{13} ((R \otimes j)U)_{23}
= \overline{U}_{13}\overline{U}_{23} \ .$$
\end{proof}

\noindent \emph{Remark}. We clearly have $\overline{\overline{U}} \simeq U.$

\vspace{2ex}

 Assume now that $(\A, \Delta)$ is  of compact type  and let $U$
    be an irreducible unitary corepresentation of $(\A_{r}, \Delta_{r}),$
    which is then necessarily finite-dimensional \cite{Wo1, Wo2}.

    \vspace{1ex}
    We will show  that the conjugate of $U$, as defined above, agrees with
the conjugate of $U$ as defined by J. Roberts
and L.\ Tuset \cite{RoTu}. We first recall their definition.

\vspace{1ex}
Let $\bar{U} \in \A \odot B(\H_U)$ be given by
$\bar{U} = (* \odot \tilde{j}) U$,
where $\tilde{j}(x) = \tilde{J} x \tilde{J}^{-1}$
for all $ x \in B(\H_U),$
and $\tilde{J}: \H_U \to \overline{\H}_U$ is any anti-linear invertible
operator such that
$\tilde{J}^*\tilde{J}$ intertwines $U$ and $(S^2 \odot \iota)U.$
 Recall here that $S^2(a) = f_{-1} * a * f_1, \ a \in \A.$
Then
$(\Delta \odot \iota)\bar{U} = \bar{U}_{13}\bar{U}_{23}$
holds as $\Delta$ is $*$-preserving and $\tilde{j}$ is multiplicative.
 The fact that $\bar{U}$ is unitary is shown in \cite{RoTu}.
(Note that $(S^2 \odot \iota)U$ is the double contragradient representation;
it is not unitary.)

    \medskip

    \begin{proposition} \label{RT}
 Assume that $(\A, \Delta)$ is  of compact type  and let $U$
    be an irreducible unitary corepresentation of $(\A_{r}, \Delta_{r}).$
   Let  $\overline{U} = (R \otimes j) U$ denote
    the conjugate of $U$ as we have defined it before. Let then
    $ \tilde{J}: \H_U \to
    \overline{\H}_U$ be
    defined by  $\tilde{J}= ( (f_{1/2} \odot j) U) J.$ Then
    $\tilde{J}$ is an anti-linear invertible operator such that
    $\tilde{J}^*\tilde{J}$ intertwines $U$ and $(S^2 \odot \iota)U.$
    Further, we have
    $\overline{U} = \bar{U},$  where $\bar{U}$ is defined as  above.
 \end{proposition}

 \begin{proof}
It is easy to check that $(f_{1/2} \odot j)(U^*)$ is invertible
    with inverse given by $ ((f_{1/2} \odot j)U^*)^{-1} = (f_{1/2} \odot
    j)U.$

    \medskip

    \noindent We set $V= I\otimes ((f_{1/2} \odot j)U^*).$ Then
    $(V^{*})^{-1} = I \otimes (f_{-1/2} \odot j) (U^*).$

    \medskip

    Recall from \cite{KuVD} that $R= S \tau_{i/2}=\tau_{i/2} S$,
    where $\tau_{i/2}(a) = f_{1/2}*a*f_{-1/2}.$

    \medskip
    \noindent Therefore we have

    \medskip

    $\overline{U} = (\tau_{i/2} \odot j) (S \odot \iota) U
    = (\tau_{i/2} \odot j)U^*$

    $$= (f_{-1/2} \odot \iota \odot f_{1/2} \odot j)
   ( (\Delta \odot \iota) \Delta) \odot \iota) U^*$$

    $$=(f_{-1/2} \odot \iota \odot f_{1/2} \odot j)
   ( (\Delta \odot \iota \odot \iota) (U_{23}^{*} U_{13}^{*})$$

   $$=(f_{-1/2} \odot \iota \odot f_{1/2} \odot j) (I \otimes I
   \otimes U^*) (\Delta \odot \iota \odot \iota) (U_{13}^*)$$

   $$= ( (f_{-1/2} \odot \iota \odot f_{1/2} \odot j)
   (\Delta \odot \iota \odot \iota) (U_{13}^*))  V$$

   $$= ( (f_{-1/2} \odot \iota \odot j)(\Delta \odot \iota)U^*)  V$$

   $$= ( (f_{-1/2} \odot \iota \odot j)( U_{23}^*  U_{13}^*))  V$$

   $$= ( (f_{-1/2} \odot \iota \odot j)( U_{13}^*)
   (f_{-1/2} \odot \iota \odot j)( U_{23}^*) )  V$$

   $$= (I \otimes (f_{-1/2} \odot j) (U^*) ) (\iota \odot j)(U^*)
   V$$

   $$= (V^*)^{-1}  (\iota \odot j)( U^*) V.$$

  From this, we get

  $$ \bar{U} = (* \odot \tilde{j})U = (* \odot  \tilde{J} \cdot
  \tilde{J}^{-1}) U = \overline{U}, $$
  (using that $(f_{1/2} \odot j) U= ((f_{-1/2} \odot j) U)^{-1} =
  ((f_{1/2} \odot j) (U^*))^{*-1},$ the first equality here relying on
  $$I_{\H_{U}}= (\varepsilon \odot \iota )U = (f_{1/2} \odot f_{-1/2}
  \odot \iota) (\Delta \odot \iota) U = (f_{1/2} \odot f_{-1/2}
  \odot \iota) (U_{13}U_{23}) $$
  $= (( f_{1/2} \odot \iota) U) ((f_{-1/2} \odot \iota) U).)$

  \medskip
  \noindent Finally, we check that $\tilde{J}^*\tilde{J}$ intertwines
  $U$ and $(S^{2}\odot \iota)U.$
  Observe first that
$$ \tilde{J}^*\tilde{J} = J^{*}((f_{1/2} \odot j)U)^{*}((f_{1/2} \odot
  j)U)J = J^{*}((f_{-1/2} \odot j)U^{*})((f_{1/2} \odot
  j)U)J $$
$$ = J^{*}((f_{-1/2} \odot j)U)^{-1}((f_{1/2} \odot
  j)U)J = J^{*}((f_{1/2} \odot j)U)((f_{1/2} \odot
  j)U)J $$
$$= J^{*}((f_{1} \odot j)U) J  = ( (f_{1} \odot \iota) U)^{*} =
F_{U}^{*}
  = F_{U}, $$
where $F_{U} = (f_{1} \odot \iota)U.$
  
\vspace{1ex} \noindent
  Now, inserting  $a = (\iota \odot \omega)U$ in $S^2(a)= f_{-1} * a *
  f_1$, we get

$$(S^2 \odot \omega)U
= S^2((\iota \odot \omega)U)
=  f_{-1} * (\iota \odot \omega)U * f_1$$
$$= ( f_1 \odot \iota \odot f_{-1} \odot \omega)
[[(\Delta \odot \iota)\Delta] \odot \iota]U$$
$$= ( f_1 \odot \iota \odot f_{-1} \odot \omega)
(\Delta \odot \iota \odot \iota) (\Delta \odot \iota)U,$$
for all $\omega \in B(\H_U)^*.$
Therefore, we have
$$(S^2 \odot \iota)U
= ( f_1 \odot \iota \odot f_{-1} \odot \iota)
(\Delta \odot \iota \odot \iota) (U_{13}U_{23})
= [I \otimes (f_1 \odot \iota)U]U [I \otimes (f_{-1} \odot
\iota)U],$$
hence
$(S^2 \odot \iota)U
= (I \otimes F_{U})U(I \otimes F_{U}^{-1}).$
Thus $\tilde{J}^*\tilde{J}=F_{U}$ intertwines $U$ and $(S^2 \odot
\iota)U,$
as claimed.

 \end{proof}

\vspace{2ex}

\subsection{Hilbert-Schmidt corepresentations}
Let $(\A,\Delta)$ be an algebraic quantum group
and $U$ be a unitary corepresentation of $(\A_r,\Delta_r).$
We introduce the Hilbert-Schmidt  corepresentation
$U_{HS}$ associated with $U$ and show that
$\overline{U} \times U \simeq U_{HS}$.

\vspace{1ex} We let $J$ and  $j$ be as in the previous subsection and denote
the
Hilbert-Schmidt operators acting on $\H_{U}$ by $HS(\H_{U}).$
We recall that  $HS(\H_U)$ is a Hilbert space  with inner product
$(x,y) = Tr (y^*x), \  x,y \in HS(\H_U),$ where $Tr $ denotes the
canonical trace on $B(\H_{U}).$

\vspace{2ex}
We define first a unitary
$\V: \overline{\H}_U \otimes \H_U \to HS(\H_U)$
by
$$\V(\eta \otimes \xi)(\xi') = (\xi', J^* \eta)_{\H_U} \, \xi
= (\eta, J\xi')_{\overline{\H}_U} \, \xi ,
\  \xi,\xi' \in \H_U, \eta \in \overline{\H}_U \ .$$

\vspace{1ex}
Define then
a normal unital $*$-isomorphism
$\tilde\V: B(\overline\H_U \otimes \H_U) \to B(HS(\H_U))$
by
$$\tilde\V(X) = \V X \V^*,
\  X \in B(\overline\H_U \otimes \H_U).$$
Note that
$\tilde\V(B_0(\overline\H_U \otimes \H_U))
= B_0(HS(\H_U)) .$
Further, let
$$\iota \otimes \tilde\V:
M(\A_r \otimes B_0(\overline\H_U \otimes \H_U))
\to M(\A_r \otimes B_0(HS(\H_U)))$$
denote the canonical extension of
$$\iota \otimes \tilde\V:
\A_r \otimes B_0(\overline\H_U \otimes \H_U)
\to \A_r \otimes B_0(HS(\H_U))
\subset M(\A_r \otimes B_0(HS(\H_U))).$$
It is clear  that $\iota \otimes \tilde\V$ is a unital $*$-isomorphism.

\vspace{1ex}
Define then
$U_{HS} \in M(\A_r \otimes B_0(HS(\H_U)))$ by
$$U_{HS} = (\iota \otimes \tilde\V)(\overline{U} \times U).$$

\begin{proposition} $U_{HS}$ is a unitary corepresentation of
$(\A_r,\Delta_r),$
with $\H_{U_{HS}} = HS(\H_U)$, which is equivalent to $\overline{U} \times
U.$
\end{proposition}

\begin{proof} $U_{HS}$ is unitary as $\iota \otimes \tilde\V$ is a unital
$*$-isomorphism and $\overline{U} \times U$ is a unitary.
Moreover,
$$(\Delta_r \otimes \iota)U_{HS}
= (\iota \otimes \iota \otimes \tilde\V) (\Delta_r \otimes \iota)
(\overline{U} \times U) = (U_{HS})_{13}(U_{HS})_{23}$$
since $\overline{U} \times U$ satisfies the corepresentation
property and $\iota \otimes \iota \otimes \tilde\V$ is multiplicative.
Finally, as $(\omega \bar\otimes \iota) U_{HS} \V
= \V (\omega \bar\otimes \iota)(\overline{U} \times U),
\  \omega \in \M_*$,
we see that $U_{HS}$ is equivalent to $\overline{U} \times U$
with unitary $\V \in Mor(U_{HS},\overline{U} \times U).$
\end{proof}

It will be useful for us later to have another way of looking at $U_{HS}.$

\vspace{1ex} \noindent
Let $l$  (resp. $r$) $: B(\H_{U}) \to B(HS(\H_{U}))$ be the normal $*$-
homomorphism (resp.\ $*$-antihomomorphism) defined by
$$l(x)(y) = xy \ \textnormal{(resp. } \ r(x)(y) = yx), \ x \in
B(\H_{U}), y \in
HS(\H_{U}).$$ It is then straightforward to check that
$$\V (I \otimes x) \V^{*} = l(x) \, , \ x \in B(\H_{U}),$$
$$\V(z \otimes I) \V^{*} = r(j^{-1}(z)) \, , \ z \in B(\overline\H_{U}).$$
Using these relations, one easily gets
$$ (I \otimes \V) X_{13} (I \otimes \V^{*}) = (\iota \bar\otimes l) X
\, \ , \ X \in B(\H) \bar \otimes B(\H_{U}),$$
$$ (I \otimes \V) Z_{12} (I \otimes \V^{*}) = (\iota \bar\otimes
rj^{-1}) Z \, \ , \ Z \in B(\H) \bar \otimes B(\overline\H_{U}).$$
Now, regarding $U \in \M \bar\otimes B(\H_{U}) \subset B(\H) \bar \otimes
B(\H_{U}),$ we have :

\begin{proposition}\label{UHS}
    $U_{HS} = (\iota \bar\otimes l) U (R \bar\otimes r) U .$
    \end{proposition}

    \begin{proof} Indeed,
 $$U_{HS}= (I \otimes \V) U_{13} \overline{U}_{12} (I \otimes \V^{*})$$
 $$= (I \otimes \V) U_{13}(I \otimes \V^{*})(I \otimes \V) \overline{U}_{12}
(I \otimes \V^{*})$$
 $$= (\iota \bar\otimes l) U  (\iota \bar\otimes rj^{-1})(R \bar\otimes j) U
 =(\iota \bar\otimes l) U (R \bar\otimes r) U .$$
 \end{proof}

\noindent \emph{Remark}. One may also associate with $U$ another Hilbert-Schmidt 
corepresentation $U_{HS'}$ of $(\A_r,\Delta_r)$ on $HS(\H_{U})$  which is given by 
$$U_{HS'} = (\iota \otimes \tilde\V \chi)( U \times  \overline{U}),$$
where $\chi$ denotes the flip map from $B(\H_{U} \otimes  
\overline\H_{U})$ to  $B( \overline\H_{U}  \otimes  
\H_{U}).$ One easily checks that  
 $U_{HS'} \simeq U \times \overline{U} \simeq 
(\overline{U})_{HS},$ and that $U_{HS'} =  (R \bar\otimes r) U (\iota 
\bar\otimes l) U.$ The two Hilbert-Schmidt corepresentations 
associated with $U$ agree 
when $\A$ is commutative, but it is unclear whether or when they are  
equivalent in the non-commutative case.

\section{Co-amenable  unitary corepresentations}
 
Inspired by \cite[Theorem 2.5]{BMT1} and \cite[Theorem 4.2]{BMT2}, we
introduce the notion of co-amenability for unitary corepresentations
of $(\A_r,\Delta_r).$

\begin{definition}
Let $(\A,\Delta)$ be an algebraic quantum group.
    A unitary corepresentation $U$ of $(\A_r,\Delta_r)$
 is said to be co-amenable if there exists $\phi \in S(\A_r)$
such that $(\phi \otimes \iota)U = I_{\H_U}$.
\end{definition}

\noindent Note that
we can equivalently require that $\phi \in S(B(\H))$
in this definition.
 The following result shows that this definition is consistent
with the notion of co-amenability for algebraic quantum groups.
Recall from \cite{BMT2,BMT3} (see \cite{BMT1} for the compact case) that an
algebraic quantum group $(\A,\Delta)$
is \emph{co-amenable} if its co-unit  $\varepsilon$ is bounded with respect
to the norm on $\A$ given by $ \| a \| = \| \pi_{r}(a) \|, a \in \A.$
Equivalently, $(\A,\Delta)$
is co-amenable if there exists a bounded linear functional
$\varepsilon_r: \A_r \to {\Complessi }$ such that
$(\iota \otimes \varepsilon_r) \Delta_r =
(\varepsilon_r \otimes \iota) \Delta_r = \iota.$ The map
$\varepsilon_r$ is then a $*$-homomorphism from $\A_{r}$ onto
$\Complessi$ and the existence of such a
 homomorphism characterizes the co-amenability of $(\A,\Delta).$

\begin{theorem} \label{allcoamen}
Let $(\A,\Delta)$ be an algebraic quantum group.
Then the following conditions are equivalent:
\begin{itemize}
\item[(1)] $(\A,\Delta)$ is co-amenable;
\item[(2)]
 $W$ is co-amenable (as a corepresentation);
\item[(3)] all unitary corepresentations
 of $(\A_r,\Delta_r)$ are co-amenable.
\end{itemize}
\end{theorem}

\begin{proof}
The equivalence (1) $\Leftrightarrow$ (2)
follows by \cite[Theorem 4.2]{BMT2}.
The implication (3) $\Rightarrow$ (2) is obvious.
In order to show (1) $\Rightarrow$ (3), we set $\phi = \varepsilon_r$
and let $U$ be a unitary corepresentations of $(\A_r,\Delta_r)$.
Then
\begin{align*}
U
& = (\iota \otimes \varepsilon_r \otimes \iota) (\Delta_r \otimes \iota) U
\\
& = (\iota \otimes \varepsilon_r \otimes \iota)(U_{13}U_{23}) \\
& = U (I \otimes (\varepsilon_r \otimes \iota)U)
\end{align*}
(here $I$ denotes the unit of $M(\A_r)$).
Multiplying by $U^*$ from the left, we get
$I \otimes I_{\H_U} = I \otimes (\varepsilon_r \otimes \iota)U$
and therefore
$(\varepsilon_r \otimes \iota)U = I_{\H_U}$.
\end{proof}

The next result may be seen as an analog of Day's classical
characterization of the amenability of a group.

\begin{proposition}
\label{coamcorep}
 let $U$ be a unitary corepresentation of $(\A_r,\Delta_r).$
Then the following conditions are equivalent:
\begin{itemize}
\item[(1)]
$U$ is co-amenable;
\item[(2)]
there exists a net $(v_i)$ of unit vectors in $\H$ such that
$$
\lim_i \ \|U(v_i \otimes \xi) - v_i \otimes \xi\|_2 = 0,
\quad \forall  \ \xi \in \H_U
$$
\end{itemize}
\end{proposition}

\begin{proof}
(2) $\Rightarrow$ (1):
By weak$^*$-compactness of $S(B(\H)),$
the net of vector states $(\omega_{v_i})$ has an accumulation point $\phi$
in $S(B(\H)).$
Passing to a subnet of $(v_i)$ if necessary, we may suppose that
$\phi(x)= \lim_i (xv_i,v_i), \ x \in B(\H)$.
 
Now, by assumption, we have
$\lim_i \ \|U(v_i \otimes \xi) - v_i \otimes \xi\|_2 = 0$,
for all $\xi \in \H_U.$
Thus
$$\omega_{\xi}((\phi \otimes \iota)U)
 = \phi((\iota \otimes \omega_{\xi})U)
= \lim_i ((\iota \otimes \omega_{\xi})(U)v_i,v_i)$$
$$= \lim_i (U(v_i \otimes \xi), v_i \otimes \xi)
= \lim_i (v_i \otimes \xi, v_i \otimes \xi)
=\omega_{\xi}(I)$$ for every $\xi \in \H_U$.
Since the set of vector states $\omega_{\xi}$ separates the elements of
$B(\H),$
it follows that $(\phi \otimes \iota)U=I$.

\vspace{1ex}
(1) $\Rightarrow$ (2):
 Let $\phi \in S(B(\H))$ be such that $(\phi \otimes \iota) U = I$.
 As $\M$ acts standardly on $\H,$ there exists a net $(v_i)$ of unit vectors
in $\H$ such that
$\phi(x) = \lim_i (xv_i,v_i), \ x \in \M$.
Then, for all $v \in \H_U$,
\begin{align*}
\lim_i (U(v_i \otimes \xi), v_i \otimes \xi)
& = \lim_i ((\iota \otimes \omega_{\xi})(U)v_i,v_i)
= \phi((\iota \otimes \omega_{\xi})U)
= \omega_{\xi} ((\phi \otimes \iota)U) \\
& = \omega_{\xi}(I)
= (\xi,\xi)
= \lim_i(v_i \otimes \xi, v_i \otimes \xi) \ .
\end{align*}
The conclusion follows easily from this.
\end{proof}

We may use the results in section 3.1 to transfer the notion of
co-amenability from corepresentations to representations : the
$*$-representation $\hat \pi_{U}$ of $\hat \A_{u}$ associated to
 a unitary corepresentation $U$ of $(\A_r,\Delta_r)$ is said to be
 \emph{co-amenable} if $U$ is co-amenable. We have for the moment  no
 intrinsic characterization of this notion.

 \vspace{1ex}
 Finally, concerning compact matrix pseudogroups \cite{Wo1}, we mention:
\begin{proposition}
Suppose that $(\A_{r},\Delta_{r})$ is a compact matrix pseudogroup
with fundamental unitary corepresentation
$U$ (so $U$ is finite-dimensional with matrix elements generating 
$\pi_{r}(\A)$ as a $*$-algebra).
Then $U$ is co-amenable
if, and only if, $(\A,\Delta)$ is co-amenable.
\end{proposition}

\begin{proof}This result is merely a restatement of \cite[Theorem
2.5]{BMT1}. A sketch of the argument is as follows.
Write $U = \sum_{i,j} u_{ij} \otimes e_{ij},$
where $u_{ij} \in \A$ and the $e_{ij}$'s form a usual system of matrix units
for $B(\H_{U}).$
If $\phi \in S(\A_r),$
then $(\phi \otimes \iota)U = I_{H_U}$
if, and only if,  $\phi(u_{ij})=\delta_{ij}$ for all $i,j$ if, and
only if, $\phi|_\A = \varepsilon.$
\end{proof}

In fact, more generally, \cite[Theorem 2.5]{BMT1} may be restated as saying
that if $(\A, \Delta)$ is of compact type and $U$ is a  unitary corepresentation of $(\A_r,\Delta_r)$ such that
its
matrix elements generate $\A_{r}$ as a C*-algebra, then $(\A,\Delta)$
is co-amenable if, and only if, $U$ is co-amenable.

 \section{Amenable unitary  corepresentations}

We first recall the following definition due to M.\  Bekka \cite{Bek}.
A continuous unitary representation $u$ of a
locally compact group $G$ on a Hilbert space $\H_{u}$ is called
\emph{amenable}
if there exists an invariant ``mean'' on $B(\H_u)$,
that is, if there exists
$m_u \in  S(B(\H_u))$ such that

$$m_u(u_{g} x u_{g}^*) = m_u(x),
\ \forall x \in  B(\H_u), \forall g \in G.
$$

This definition has no obvious counterpart in the quantum group setting.
However, Bekka also introduces a notion of ``topological'' invariant
mean whose existence is equivalent to the amenability of $u$,
see \cite[Theorem 3.5]{Bek}.
Inspired by this result, we introduce the following notion:

\begin{definition}
A unitary corepresentation $U$
of $(\A_r,\Delta_r)$
is called left-amenable (resp.\ right-amenable) if there exists
$m_U$ (resp.\ $m'_{U})$) $\in S(B(\H_U))$ such that

$$m_U((\omega \bar\otimes \iota)(U^* (I \otimes x) U)) = \omega(I)
m_U(x)$$

$$(\textnormal{resp.} \ m'_U((\omega \bar\otimes \iota)(U (I \otimes
x) U^*)) = \omega(I)
m'_U(x))$$
for all $ x \in B(\H_U), \omega \in \M_{*}.$

The state $m_U$ (resp. $m'_{U}$) is called a left-invariant (resp.\
right-invariant) mean for $U$.
\end{definition}

 \emph{Remarks.}

\smallskip
$(i)$ When $U$ is the unitary corepresentation of
$(C_{0}(\Gamma),\Delta)$  associated to a unitary representation $u$ of a
discrete group $\Gamma,$ one easily checks that $U$ is left-amenable
(resp. $U$ is right-amenable)
if,
and only if, $u$ is amenable. This is a simple consequence of the fact
that $(\delta_{\gamma} \bar\otimes \iota) U = u_{\gamma},$ where the
delta function at $\gamma \in \Gamma,$
$\delta_{\gamma},$ is considered as an element of $\ell^1(\Gamma),$
that is, of the predual of
$\ell^{\infty}(\Gamma).$ Further, in this case, it is quite obvious that  a
left-
(resp.\ right-) invariant mean for
$U$ is both left- and right-invariant.

$(ii)$ We don't know
whether the existence of a left-invariant mean
for $U$ is equivalent to the existence of a right-invariant one in the
general situation. However, we have
 $U$ is left-(resp.\ right-) amenable if, and only if, $\overline{U}$ is
 right-(resp.\ left-) amenable.

 Indeed, if  $m_U$ is a left-invariant mean for $U,$
then
$ m_{U} \circ j^{-1}$ is a right-invariant mean
for $\overline
U.$ If $m'_{U}$ is a right-invariant mean
for $\overline
U,$ then $m'_{U} \circ j$ is a left-invariant mean for $U.$
The resp.\ assertions are proven similarly.

\smallskip
$(iii)$ The property of left-amenability (resp. right-amenability) is clearly invariant under unitary
equivalence.

\smallskip
$(iv)$  By ``linearizing'' the concept of amenability, one gets
a  related, but seemingly independent, notion:
a unitary corepresentation $U$ of $(\A_r,\Delta_r)$
is said to be \emph{hypertracial} if
there exists $m''_U \in S(B(\H_U))$ such that
\begin{equation}
\label{hypertracecorep}
m''_U((\omega \bar\otimes \iota)(U (I \otimes x)))
= m''_U((\omega \bar\otimes \iota)((I \otimes x)U)),
\ \forall x \in B(\H_U), \forall \omega \in \M_{*}.
\end{equation}
Actually, condition (\ref{hypertracecorep}) is equivalent to
$$m''_U((\omega \bar\otimes \iota)(U) x) =
m''_U(x (\omega \bar\otimes \iota)(U))
\ \forall x \in B(\H_U), \forall \omega \in \M_{*},$$
which in turn is equivalent to
$$m'_U(\hat \pi_U(a)x) = m'_U(x \hat \pi_U(a)),
\forall x \in B(\H_U), \ \forall a \in \A_u \ .$$
Hence, hypertraciality of $U$ is equivalent to hypertraciality of
$\hat \pi_{U}$ in the sense of \cite{Bed}.
This hypertrace property is easily seen to
correspond  to left- and right-amenability in the  case of a
corepresentation
arising from a unitary
representation of a discrete group.
 
\vspace{2ex}
Now recall from \cite{BMT2,BMT3} that an algebraic quantum group
$(\A,\Delta)$ is called  \emph{amenable} if
there exists a left-invariant mean for $(\A,\Delta),$
that is, if there exists  $m \in S(\M)$ such that
$$m((\omega \bar\otimes \iota)\Delta_{r}(x)) = \omega(I) m(x), \
\forall x \in \M, \forall \omega \in \M_{*}.$$
A right-invariant mean for $(\A,\Delta)$ is defined similarly.
By composing with the anti-unitary antipode $R$
(which is defined on $\M$, see \cite{KuVD,KV}),
one easily sees that the existence of a
left-invariant mean for $(\A,\Delta)$ is equivalent to
the existence of a right-invariant mean for it.
It is then  straightforward to check that $(\A,\Delta_{\rop})$ is   amenable
if and only if $(\A,\Delta)$ is amenable.

\medskip
The following result is well known (see \cite{Bek} for the
equivalence between (3), (4) and (5); the equivalence between
(1), (2) and (3) is merely classical, as explained in \cite{BMT1,BMT2}).

\begin{theorem}
Let $\Gamma$ be a discrete group and let $(\A,\Delta)$ be the algebraic
quantum group associated with $\A = {\Bbb C}[\Gamma],$ the
group-algebra of $\Gamma,$ so $\hat \A = C_{c}(\Gamma).$
Then the following are equivalent:
\begin{itemize}
\item[(1)] $(\hat\A,\hat\Delta)$ is amenable;
\item[(2)] $(\A,\Delta)$ is co-amenable;
\item[(3)] $\Gamma$ is amenable;
\item[(4)] the (left-) regular representation $\lambda$ of $\Gamma$ is
amenable;
\item[(5)] all unitary representations of $\Gamma$ are amenable.
\end{itemize}
\end{theorem}

In the case of an algebraic quantum group, it is known \cite{BMT2} that
$(1)$ implies
$(2).$ The converse implication is only known to hold when
$(\A,\Delta)$ is compact with a tracial Haar functional \cite{BMT3}.
We will give another proof of this result in Section 9.
However, amenability of an algebraic quantum group may be characterized
through amenability of its corepresentations as follows.

\begin{theorem} \label{allamen}
Let $(\A,\Delta)$ be an algebraic quantum group.
Then the following conditions are equivalent:
\begin{itemize}
\item[(1)] $(\A,\Delta)$ is amenable;
\item[(2)] $W$
 is left-amenable (as a corepresentation);
\item[(3)] $\overline W$ is right-amenable (as a corepresentation);
\item[(4)] all unitary corepresentations of $(\A_r,\Delta_r)$
are left-amenable.
\item[(5)] all unitary corepresentations of $(\A_r,\Delta_r)$
are right-amenable.
\end{itemize}
\end{theorem}

\begin{proof}
The implications (4) $\Rightarrow$ (2) and (5) $\Rightarrow$ (3) are
obvious.
The equivalence between (2) and (3) is just a special case of $(ii)$ in
our previous remark.

\vspace{1ex}
(2) $\Rightarrow$ (1):
Let $m_W \in S(B(\H))$ be a left-invariant mean for $W,$
and let $m$ to be the restriction of $m_W$ to $\M$.
Then $m$ is clearly a state, which is left-invariant since
$$m((\omega \bar\otimes \iota)\Delta_r(x))
= m_W((\omega \bar\otimes \iota)(W^* (I \otimes x) W))
= \omega(I) m_W(x) = \omega(I) m(x),$$
for all $x \in \M, \omega \in \M_{*} \ .$

\vspace{1ex}
(1) $\Rightarrow$ (4) and (1) $\Rightarrow$ (5):
Let $m$ be a left-invariant mean for $(\A,\Delta)$
and let $U$ be a unitary corepresentation of  $(\A_r,\Delta_r).$
We pick a normal state $\Omega$ on $B(\H_U)$ and define
$m_U \in S(B(\H_U))$ by
$$m_U(x)= m ((\iota \bar\otimes \Omega)(U^* (I \otimes x) U)), \ \ x
\in B(\H_{U}).$$
Then we check the validity of the equation expressing the
left-invariance    property
of $m_{U}:$
the {\it l.h.s.} is
$$m_{U}(x) =m(\iota\bar\otimes \Omega)
(U^*(I \otimes [(\omega \bar\otimes \iota)(U^*(I \otimes x)U)])U)$$
while
  the {\it r.h.s.} is equal to
$$\omega(I)m_{U}(x)=\omega(I)m ((\iota \bar\otimes \Omega)(U^* (I \otimes
x) U))$$
$$= m(\omega \bar\otimes \iota)\Delta((\iota \bar\otimes \Omega)
(U^* (I \otimes x)U)) \ (\textnormal{using left-invariance of}\ m)$$
$$ = m(\omega \bar\otimes \iota)(\iota \bar\otimes \iota \bar\otimes \Omega)
(\Delta \bar\otimes \iota) (U^* (I \otimes x)U) $$
$$ =  m(\omega \bar\otimes \iota)(\iota \bar\otimes \iota \bar\otimes
\Omega)
(((\Delta \bar\otimes \iota)U)^* (I \otimes I \otimes x)
(\Delta \bar\otimes \iota)U) $$
$$=  m(\omega \bar\otimes \iota)(\iota \bar\otimes \iota \bar\otimes \Omega)
((U_{13}U_{23})^* (I \otimes I \otimes x) U_{13}U_{23}) $$
$$ = m(\iota \bar\otimes \Omega) (\omega  \bar\otimes \iota \bar\otimes
\iota)
(U^*_{23}U^*_{13}(I \otimes I \otimes x)U_{13}U_{23}).$$

Hence, the desired conclusion follows from the identity
$$
U^*(I \otimes [(\omega \bar\otimes \iota) (U^* (I \otimes x) U)]  )U
=
(\omega \bar\otimes \iota \bar\otimes \iota)
(U_{23}^*U_{13}^* (I \otimes I \otimes x) U_{13}U_{23})
$$
which is  easily verified.

Similarly, we define $m'_{U} \in S(B(\H_U))$ by

$$m'_U(x)= m\circ R ((\iota \bar\otimes \Omega)(U (I \otimes x)
U^*)),  \ \ x
\in B(\H_{U}).$$
Then, using  the right-invariance of $m \circ R,$ one now checks that
$m'_{U}$ is a right-invariant mean for $U.$

\vspace{1ex}
\noindent It clearly follows that all stated conditions  are equivalent.

\end{proof}

We know from \cite[Theorem 4.7]{BMT2} that co-amenability of
$(\A,\Delta)$ implies that $(\hat \A, \hat \Delta)$ is amenable,
hence that $(\hat \A, \hat \Delta_{op})$ is amenable.
By combining this fact with Theorem \ref{allcoamen} and Theorem
\ref{allamen},
we see that if all the unitary corepresentations of $(\A_r,\Delta_r)$
are co-amenable then all the unitary corepresentations of
$(\hat\A_r,\hat\Delta_{r,\rop})$ are amenable.
This lends some evidence  that there might be
some
 correspondence between co-amenable elements in
Corep$(\A_r,\Delta_r)$ and amenable elements in
Corep$(\hat\A_r,\hat\Delta_{r,\rop}).$

\vspace{2ex}

By using Theorem \ref{cateq}, one may clearly transfer the notion of
amenability to representations of algebraic quantum groups. Theorem
\ref{allamen} may then be reformulated in an obvious manner.

\vspace{2ex}
 
An analog of \cite[Theorem 3.6]{Bek}, which characterizes the amenability
of a unitary representation of a group, is as follows.

\begin{proposition}\label{BekTC}
    Let $U$ be a unitary corepresentation of
$(\A_{r},\Delta_{r}).$ Organize $TC(\H_{U})$, the trace class
operators on $\H_{U},$ as a Banach $\M_{*}$-module by means of

$$Tr ((\omega \cdot s)x) = Tr ( s(\omega \bar{\otimes} \iota)(U^* (I
\otimes x) U)), $$
$\omega \in \M_{*}, s \in TC(\H_{U}), x \in B(\H_{U}).$

\vspace{1ex} \noindent
Then $U$ is left-amenable if, and only if, there exists a net
$(s_{i})$ in $TC(\H_{U})^{+}_{1}$ such that
$$ \lim_{i} \| \omega \cdot s_{i} - s_{i} \|_{1} = 0,
\quad \forall \omega \in \M_{*}.$$

    \end{proposition}

    \begin{proof} The proof is an easy adaptation of
 the proof of \cite[Theorem 3.6]{Bek}. If $(s_{i})$ is a net as above,
 then a left-invariant mean for $U$ is obtained by picking any
 weak$^{*}$-limit point of the net $(m_{i}) \subset S(B(\H_{U}))$ given by
 $m_{i}(\cdot) = Tr(s_{i} \cdot).$ Conversely, assume that $m_{U}$ is
 a left-invariant mean for $U.$ As the normal states are
 weak$^{*}$-dense in $S(B(\H_{U}))$, we may pick a net $(t_{i}) \subset
 TC(\H_{U})^{+}_{1}$ such that $m_{U}$ is  weak$^{*}$-limit point of
 the net $(Tr(t_{i} \, \cdot)) \subset S(B(\H_{U})).$ Namioka's classical
 argument \cite{Pa} gives then the existence of a net $(s_{i})$ with the
 required properties.

 \end{proof}

    One may clearly also obtain a similar characterization of
    right-amenability for unitary corepresentations of $(\A_{r}, 
    \Delta_{r}).$

\vspace{2ex}
To illustrate the notion of invariant mean for corepresentations, we
now consider the case
where $(\A,\Delta)$ is of compact type. Let then $U$ be a
finite-dimensional unitary representation of $(\A_{r},\Delta_{r}).$
As $(\A,\Delta)$ is amenable, see the paragraph preceding Theorem 4.7
in \cite{BMT2}, we deduce from Theorem \ref{allamen} that all  unitary
corepresentations of $(\A_{r},\Delta_{r})$ are left- (and right-) amenable.
We shall now describe
somewhat more explicitly
 a left-invariant mean $m_U$ for $U$, following the construction given in the
proof of Theorem
\ref{allamen}.
 
 \vspace{1ex}
 Identifying $\A$ with the dense Hopf *-subalgebra
 $\pi_{r}(\A)$ of
 $(\A_{r},\Delta_{r}),$ we may write
 $U = \sum_i a_i \otimes b_i \in \A \odot B(\H_U)$ for some $a_{1},
 \ldots , a_{N} \in \A,  b_{1}, \ldots, ,b_{N} \in B(\H_{U}).$
 \vspace{1ex}
 Recall that a left-invariant mean for $U$ is provided by
 $$m_U(x)= \varphi_{r}((\iota\bar\otimes \Omega)(U^*(I \otimes x)U) ),
\ x \in B(\H),$$ where we have the freedom to choose any
$\Omega \in B(\H_U)^+_{*,1}$.
Set $d_{U}= \dim \H_U$
and let $\tau= 1/d_U {\rm Tr}, $
denote the normalized trace on $B(\H_U).$
Plugging in $\Omega = \tau,$ we get
\begin{align*}
m_U(x) & = \varphi_{r}( (\iota\bar\otimes \tau)(U^*(I \otimes x)U) ) \\
& = \sum_{i,j} \varphi (\iota \otimes \tau) (a_i^*a_j \otimes b_i^* x b_j))
\\
& = \sum \varphi(a_i^*a_j)\tau( b_i^* x b_j) \\
& = 1/d_U{\rm Tr} (\sum_{i,j}\varphi (a_i^*a_j)b_jb_i^* x)
\end{align*}
so that $m_U(\cdot) = {\rm Tr}(K_{U} \cdot),$
 where the density matrix $K_U \in B(\H_U)$ is given by
\begin{equation}
K_{U} = 1/d_U (\varphi \odot\iota)(U (\sigma \otimes \iota) U^*),
\end{equation}
where $\sigma$ is the automorphism of $\A$ given by $\sigma(a)=f_1 * a *
f_1, \ a \in \A.$
 
We remark that if $\varphi$ is tracial, then $f_{1} = \varepsilon,$
hence $\sigma = \iota$ and thereby
  $K = I/d_U,$
that is, the left-invariant mean for $U$ is just $\tau.$

\vspace{2ex}
Now assume that $U$ is  irreducible.
We write $U = \sum_{i,j} u_{ij} \otimes m_{ji}$,
where $u_{ij} \in \A$ and the $m_{ij}$'s form a system of matrix units for
$B(\H_{U})$
such that $m_{lk}m_{sr} = \delta_{lr}m_{sk}$ and $m^*_{lk}=m_{kl}$.
 Using the orthogonality relation
$$\varphi((u_{km})^*u_{ln})
= (1/M_{U}) \delta_{mn}f_{-1}(u_{lk}) \ ,$$
where $M_{U}$
 denotes  the quantum dimension of $U,$ we get
 Then
\begin{align*}
\sum_{i,j,k,l} \varphi((u_{ij})^*u_{kl})m_{lk}m^*_{ji}
& =  \sum_{i,j,k,l} (1/M_{U}) \delta_{jl} f_{-1} (u_{ki}) \delta_{lj} m_{ik}
\\
& = \sum_{i,j,k} (1/M_{U}) f_{-1} (u_{ki})m_{ik} \\
& = (d_U/M_{U}) \sum_{i,k} f_{-1} (u_{ki})m_{ik} \\
& =  (d_U/M_U) (f_{-1} \odot \iota) (\sum_{i,k} u_{ki}\otimes m_{ik}) \\
& = (d_U/M_U)(f_{-1} \odot \iota)(U) \ .
\end{align*}
Hence, in this case, we get
 $$K_U = (f_{-1} \odot \iota)(U)/M_U \ .$$

We summarize what we have shown.

\begin{proposition} \label{meanirr}
  Assume that $(\A,\Delta)$ is of compact type and let $U$ be a
finite-dimensional unitary representation of $(\A_{r},\Delta_{r}).$
Let $d_{U}$ (resp.\ $\M_{U}$ ) denote the usual (resp.\ quantum)
dimension of $U,$ and let $\sigma$ be
the automorphism of $\A$ given by $\sigma(a)=f_1 * a * f_1, \ a \in \A.$
Then a left-invariant mean $m_{U}$ for $U$ is given by $m_{U}(\cdot)=$
Tr$(K_{U} \, \cdot),$ with density matrix $K_{U}$ given by
$$K_{U} = 1/d_U (\varphi \odot \iota)(U (\sigma \otimes \iota)
U^*) . $$
If $U$ is is irreducible, then
$$K_{U} = (f_{-1} \odot \iota)(U)/M_U.$$
\end{proposition}

\section{On weak containment}
We discuss  in this section the notion of weak containment
for representations and corepresentations of  algebraic quantum groups.
We begin by discussing the stronger
(and easier)
notion of containment.

\subsection{Strong  containment}

 We recall the following definition.

\begin{definition}Let $(\A,\Delta)$ be an algebraic quantum group and let
$\pi_1, \pi_2$ be two non-degenerate
$*$-representations of $\A_u.$
We say that $\pi_1$ is \emph{contained} in $\pi_2$,
and write $\pi_1 < \pi_2$,
if there exists an isometry $T \in Mor(\pi_1,\pi_2)$.
\end{definition}

Observe that $K = T (\H_{\pi_1}) $ is then a closed invariant subspace for
$\pi_2.$
 Therefore, if $\pi_{2}$ is irreducible, then
any non-degenerate $*$-representation $\pi_{1}$ of $\A_u$ contained in
$\pi_{2}$
is unitarily equivalent to $\pi_{2}$.

The interesting case where $\pi_{1}=\varepsilon_{u}$ may be characterized as
follows.

\begin{proposition} \label{strongcont}
Let $(\A,\Delta)$ be an algebraic quantum group and consider
a non-degenerate $*$-representation $\pi$ of $A_u.$
Write
$\pi =\pi_U$ for a unique unitary corepresentation $U$
of $(\hat\A_r ,\hat\Delta_{r,\rop}).$
The following conditions are equivalent:
\begin{itemize}
\item[(1)] $\varepsilon_u < \pi_U$;
\item[(2)] there exists a unit vector $\xi \in \H_U$ such that
$(\iota \bar\otimes \omega_\xi)U = I \in M(\hat\A_r)$;
\item[(3)] there exists a unit vector $\xi$ in $H_U$
such that
$U(v \otimes \xi) = v \otimes \xi, \ \forall v \in \H \ .$
\end{itemize}
\end{proposition}

\begin{proof}
We first show $(1)$ implies $(2)$. Assume that $(1)$ holds.  Then there
exists a linear map
$T: {\Bbb C}\rightarrow\H_U$ such that $T^* T =1$ and
$T\varepsilon_u (a)=\pi_U (a)T$ for all $a\in\A_u$.
Consider the unit vector $\xi =T(1)$ of $\H_U$.
Then the adjoint $T^* :\H_U\rightarrow {\Bbb C}$
is given by $T^* (\eta )=(\eta ,\xi )$ for all $\eta\in\H_U$.
Now, for all $a\in \A\subset \A_u$, we have
\begin{align*}
Q (a)(I)
& =(Q (a)\bar{\otimes}\iota)(I\otimes 1)
=\varepsilon (a)
=(\varepsilon_u (a)1,1)_{\Bbb C} \\
& =(T^*\pi_U (a)T(1),1)_{\Bbb C}
=(\pi_U (a)T(1),T(1))_{\H_U} \\
& =(\pi_U (a)\xi ,\xi )_{\H_U}
=((Q (a)\bar{\otimes}\iota)U\xi ,\xi )_{\H_U}
= Q(a) ((\iota \bar\otimes \omega_\xi)U) \ .
\end{align*}
Since $Q(\A)$ is dense in $\hat\M_*$ it follows that
$I = (\iota \bar\otimes \omega_\xi)U$.

\vspace{1ex}
Next, we show that $(2)$ implies $(1).$ Given a unit vector $\xi$
satisfying $(2),$
we define the linear isometry $T: {\Bbb C} \to \H_U$ by
$T(1)= \xi.$
By reversing the above calculations, we see that
$\varepsilon_u(a)=T^* \pi_U(a)T$ holds for all $a \in \A_u.$
It is easily checked that $T \in Mor(\varepsilon_u,\pi_U).$

\vspace{1ex}
Finally, to prove the equivalence between (2) and (3), observe first that
$(\iota \bar\otimes \omega_\xi)U = I$ if and only if
$(\omega_v \bar\otimes \omega_\xi)U = 1$ for all unit vectors
$v \in \H.$
Now, for a unit vector $v\in\H,$ one easily checks that
$$(U(v\otimes\xi),v\otimes\xi )
=1\Leftrightarrow \|U(v\otimes\xi )-v\otimes\xi\|_2 =0.$$
Hence, this equivalence is clear, and the proof is finished.

\end{proof}

Let $U,V$ be unitary corepresentations of $(\A_{r},\Delta_{r}).$
We say that $U$ is (strongly) \emph{contained} in $V,$
 and write $U < V$ if there is an isometry $T \in  Mor(U,V).$ It is an
 easy exercise to check that $U < V$ if, and only if $\hat{\pi}_{U} <
 \hat{\pi}_{V}.$ One may then clearly obtain a  result
 similar to Proposition \ref{strongcont}.

 \vspace{1ex} \noindent
\emph{Example.} Let $(\A,\Delta)$ be of compact type and  $U$ be an
irreducible unitary representation of $(\A_{r},\Delta_{r}).$ Let $\{e_{i}\}$
denote an
orthonormal basis for $\H_{U}$ and $\tilde{J}$ be defined as in
Proposition \ref{RT}. Then the isometry $R$ from $\Complessi$
into $\overline{\H}_{U} \otimes \H_{U}$ determined by
$$R(1) = \sum_i {\tilde{J}^*}{}^{-1}(e_i) \otimes e_i \ $$
satisfies $R \in Mor ( I \otimes 1,\overline{U} \times U),$ as shown
in \cite{RoTu}. Hence, $I \otimes 1 < \overline{U} \times U.$ It follows
that
$\xi = \sum_i {\tilde{J}^*}{}^{-1}(e_i) \otimes e_i
\in \overline{\H}_U \otimes \H_U$
satisfies $(\overline{U} \times U)(\eta \otimes \xi) = \eta \otimes
\xi$ for all $\eta \in \H.$
Since $(\A, \Delta)$ is compact by assumption, we know that
all
unitary corepresentations of $(\A_{r},\Delta_{r})$ are
left- (and right-) amenable, as pointed out in the previous section.
Indeed, we have  seen in Proposition \ref{meanirr} that a left-invariant
mean $m_U$ for $U$ is given by
$$m_U(x) = Tr(K_U x), \  x \in B(\H_U),$$
where $K_U = \frac{1}{M_{U}}(f_{-1} \odot \iota)U.$
 
\vspace{1ex}
\noindent Now, let $$\widetilde{\V(\xi)}= \frac{1}{(M_{U})^{1/2}} \V(\xi),$$
where $\V : \overline{\H}_{U} \otimes \H_{U} \to HS(\H_{U})$ is
defined as in Section 4.
 It is then not difficult to check by direct  computation that we also have
$$m_{U} ( \cdot ) = Tr(\widetilde{\V(\xi)}^* \cdot
\widetilde{\V(\xi)}),$$
that is,  we have $\ \ \widetilde{\V(\xi)} \widetilde{\V(\xi)}^* = K_{U}.$

\subsection{Weak containment}

Let $(\A,\D)$ be an algebraic quantum group  and let $\pi_1, \pi_2$ be
non-degenerate $*$-representations of $\A_u.$
As usual for representations of C*-algebras,
we say that $\pi_1$ is \emph{weakly contained} in $\pi_2,$
and write $\pi_1 \prec \pi_2,$
if Ker $\pi_2 \subset$ Ker $\pi_1.$
This relation is obviously transitive and reflexive, and,
of course, $\pi_1 < \pi_2$ implies $\pi_1 \prec \pi_2$.

\begin{proposition}
 With notation as above, we have $\pi_1 \prec \pi_2$
if and only if there exists a unique surjective $*$-homomorphism
$\theta: \pi_2(\A_u) \to \pi_1(\A_u)$ such that
$\theta\pi_2(a) = \pi_1(a), \ \forall a \in \A.$
\end{proposition}

\begin{proof}
This  proof is easy and left to the reader.
\end{proof}

An almost immediate consequence of this proposition is  the following.

\begin{corollary}
\label{charcoamen1}
$\varepsilon_u \prec \pi_r$ if and only if $(\A,\Delta)$ is co-amenable.
\end{corollary}

\emph{Remark.} Let $(\A, \Delta)$ be an algebraic quantum group.
Note that
$(\A,\Delta)$ is co-amenable
if, and only if, $\pi \prec \pi_r$ for every
non-degenerate $*$-representation $\pi$ of $\A_{u}$.
Indeed, if $(\A,\Delta)$ is co-amenable, then $\A_u = \A_r,$
that is, $\pi_r: \A_u \to \A_r$ is injective (see \cite{BMT2}).
Therefore, $\{ 0 \} $ = Ker $\pi_r  \subset $ Ker $\pi.$ On
the
other hand, if
$(\A,\Delta)$ is not co-amenable, then $\varepsilon_u$ is not weakly
contained in $\pi_{r}.$

We also  remark that  the condition
$\varepsilon_u \prec \pi$ and
the condition $ \pi \prec \pi_{r}$ are generally independent of each
other. In fact, if $\pi = \varepsilon_u,$ then
the first is trivially satisfied, while the second holds if and only
if $(\A, \Delta)$ is co-amenable. On the other hand, if $\pi =
\pi_{r},$ then second is trivially satisfied, while the first holds
if and only if if $(\A, \Delta)$ is co-amenable.

\begin{definition}
\label{WCPrep}Let $(\A,\D)$ be an algebraic quantum group.
A non-degenerate $*$-representation $\pi$ of $\A_u$
 is said to have the weak containment property (WCP)
if $\varepsilon_u \prec \pi$,
that is, Ker $\pi \subset$ Ker $\varepsilon_u$.
\end{definition}

Thus $\pi$
has the WCP if and only if
there exists a $*$-homomorphism $\theta: \pi(\A_u) \to {\Bbb C}$
such that
$\theta \pi(a) = \varepsilon(a)$ for all $a \in \A \subset \A_u.$
 
\begin{definition}
\label{WCPcorep}
Let $(\A,\Delta)$ be an algebraic quantum group, $U, V$ be
unitary corepresentations
of $(\hat\A_r,\hat\Delta_{r,\rop})$ and let $ \pi_{U}, \pi_{V}$ be the
associated  $*$-representations of $(\A_{u},\Delta_u).$
 
We say that $U$ is weakly contained in $V$ if $ \pi_{U}$ is weakly
contained in $\pi_{V}.$ Moreover, we say that $U$
has the weak containment property (WCP) if the trivial
corepresentation $I \otimes 1$ is weakly contained in $ U,$ that is,
if $ \pi_{U}$ has the WCP.
\end{definition}

\begin{corollary}
\label{charcoamen2}
An algebraic quantum group $(\A,\Delta)$ is co-amenable
if and only if
$\hat{W},$ as a unitary corepresentation of
$(\hat\A_r,\hat\Delta_{r,\rop}),$
has the WCP.
\end{corollary}

\begin{proof} As $\pi_{\hat{W}} = \pi_{r}$, see Theorem \ref{cateq},
this is just a reformulation of Corollary \ref{charcoamen1}.
\end{proof}

The weak containment property for unitary corepresentations may be
characterized as follows.

\begin{theorem}
\label{WCPequiv}
Let $(\A,\Delta)$ be an algebraic quantum group, $U$ be
a unitary corepresentation
of $(\hat\A_r,\hat\Delta_{r,\rop})$ and let $ \pi_{U}$ be the
associated  $*$-representation of $\A_{u}.$
 The following conditions are equivalent:
\begin{itemize}
\item[(1)] $I \otimes 1 \prec U,$ that is, $U$ has the WCP ;
 \item[(2)] there exists $\psi \in S(B(\H_U))$ such that
$\psi((\omega\bar{\otimes}\iota )U)=\omega(I), \ \forall
\omega\in\hat\M_{*}$;
\item[(3)] there exists a net $(\xi_i)$ of unit vectors in $H_U$
such that
$$
\lim_i \ \|U(v \otimes \xi_i) - v \otimes \xi_i\|_2 = 0 \ \forall v \in \H \
;
$$
\item[(4)] there exists a net $(\xi_i)$ of unit vectors in $H_U$
such that
$$
\lim_i \ (U(v \otimes \xi_i),v \otimes \xi_i) = 1, \ \forall v \in \H,
\|v\|_2=1 \ .
$$
\end{itemize}
Further, any of these conditions implies that $U$ is left-amenable,
right-amenable
and hypertracial.
\end{theorem}

\begin{proof} The equivalence between (3) and (4) is elementary.

    \vspace{1ex}
(1) $\Rightarrow$ (2) and (1) $\Rightarrow$ (4):
Assume that $I \otimes 1 \prec U.$
By the remark following Definition \ref{WCPrep}
there exists a $*$-homomorphism
$\theta: \pi_U(\A_u) \to {\Bbb C}$ such that
$\theta(\pi_U(x)) = \varepsilon_u(x)$ for all $x \in \A_u$.
We extend the state $\theta$ to a state $\psi$ on $B(\H_U)$.
Then, for all $a \in \A \subset \A_u$, we have
$$\psi((Q(a) \bar\otimes \iota)U)=\psi((Q_{r}(a) \otimes \iota)U)
= \psi(\pi_U(a))$$
$$=  \varepsilon_u(a)
= (Q_r(a) \otimes \iota)(I \otimes 1)
= Q_r(a)(I)= Q(a)( I).$$
Since $Q(\A)$ is dense in $\hat\M_*,$
 we  get, by continuity,
$\psi((\omega\bar{\otimes}\iota )U)=\omega(I),
\ \forall \omega\in\hat\M_{*}$,
which shows (2). \\

Further, as $\varepsilon_u$ is a $*$-homomorphism on $\A_{u},$ it is a
pure state on $\A_{u}.$ From \cite[Proposition 3.4.2, ii)]{Dix}, we can
then conclude that there exists a  net of unit vectors $(\xi_{i}) \in
\H_{U}$ such that $ \varepsilon_u (x) = \lim_{i} ( \pi_{U}(x)
\xi_{i}, \xi_{i})$ for all $ x \in \A_{u}.$ Since  $ \varepsilon_u =
\psi \circ \pi_{U}$ as above, this means that $\psi (y) = \lim_{i}
\omega_{\xi_{i}}(y) $ for all $y \in \pi_{U}(\A_{u}).$ As $(\omega
\bar \otimes \iota) U \in   \pi_{U}(\A_{u})$ (see \cite[Theorem 3.3]{BMT3})
and
$\psi((\omega\bar{\otimes}\iota )U)= \omega(I),$ for all
$\omega\in\hat\M_{*},$
we get $\lim_{i} \omega_{\xi_{i}} (( \omega_{\eta} \bar \otimes \iota)
U) = 1 $ for all unit vectors $\eta$ in $\H.$
This  just says that $ \lim_{i} (U (\eta \otimes \xi_{i}) , \eta \otimes
\xi_{i})=1$
for all unit vectors $\eta$ in $\H,$
hence that (4) holds.

\vspace{1ex}
(2) $\Rightarrow$ (1): Assume that (2) holds, and let $\psi$ be as
in (2).
Let  $x \in $ Ker $\pi_U.$ Choose a sequence $(a_n)$ in $\A$
converging to $x \in \A_u$ with respect to the norm $\|\cdot\|_u$.
Then, by continuity of $\pi_U$, we get
$(Q_r(a_n) \otimes \iota)U = \pi_U(a_n) \to \pi_U(x) = 0$.
Using the assumption,  we have
$\psi((Q(a_n) \bar\otimes \iota)U)=Q(a_n)(I)$ for all $n.$
By continuity of $\psi$ we therefore get
$$\varepsilon_u(a_n) = Q_r(a_n)(I)= Q(a_n)(I)$$
$$= \psi((Q(a_n) \bar\otimes \iota)U)= \psi((Q_{r}(a_n) \otimes
\iota)U) \to 0.$$
Thus, by continuity of $\varepsilon_u$,
we get $\varepsilon_u(x) = \lim_n \varepsilon_u(a_n) = 0$,
so $x \in $ Ker $\varepsilon_u$. Hence, (1) holds.

\vspace{1ex}
(4) $\Rightarrow$ (2): Let  $(\xi_i )_i$ be a net satisfying condition (4).
 Using Alaoglu's theorem, and passing to a subnet if necessary, there
exists a
$\psi\in S(B(H_U ))$ such that $\psi (x)=\lim_i \omega_{\xi_i} (x)$, for all
$x\in B(H_U )$. Since $\hat M$ is in  standard form on $H$, any normal
state $\omega$ on $\hat M$ is of the form $\omega_v$ for some unit vector
$v\in \H.$
Then $\psi((\omega_v \bar\otimes \iota)U)
=\lim_i \omega_{\xi_i} ((\omega_v \bar\otimes \iota)U)
=\lim_i \ (U(v \otimes \xi_i),v \otimes \xi_i) = 1 = \omega_{v}(I),$
so $\psi$ satisfies condition (2).

\vspace{1ex}
\noindent Hence, we have established the equivalence between
conditions (1)-(4).

\vspace{1ex}
Finally, assume that (2) holds and
set $m_U = \psi.$ Let $\omega \in \hat{\M}^{+}_{*,1}.$ Then
$m_U ((\omega \bar\otimes \iota) U)
 = \omega(I) = 1.$
As $U$ is a unitary in $\hat{\M} \bar\otimes B(\H_{U}),$ it follows
from the Cauchy-Schwarz inequality  that the
state
$m_U ((\omega \bar\otimes \iota) (\cdot))$  on
$\hat \M \bar\otimes B(\H_{U})$ is multiplicative at $U$ and at $U^{*}.$
Hence,
$$m_U((\omega \bar\otimes \iota)(U^*(I \otimes x)U))
= m_U ((\omega \bar\otimes \iota) U^*)
m_U ((\omega \bar\otimes \iota) (I \otimes x))
m_U ((\omega \bar\otimes \iota) U)$$
$$= m_U(x) = \omega(I) m_U(x) $$
for all $x \in  B(\H_U)$ and $\omega \in \hat{\M}^+_{*,1}.$
It easily follows that  $m_{U}$ is a left-invariant mean for $U.$ Similarly,
$m_{U}$ is  a right-invariant mean for $U,$ and it also serves to
show that $U$ is hypertracial. This finishes the proof.
\end{proof}

Weak containement and WCP for unitary corepresentations of $(\A_{r},
\D_{r})$
are defined in an analogous way, via weak containement and WCP for the
associated
representations of $\hat{\A}_{u}.$ From a conceptual point of view,
it is   better to work in this setting, and we will often do this in
the sequel.  All statements concerning
WCP for
unitary corepresentations of $(\hat\A_r,\hat\Delta_{r,\rop})$ have an
analogous statement concerning WCP for
unitary corepresentations of $(\A_{r}, \D_{r}).$ For example, we have
the following counterpart to Theorem \ref{WCPequiv}.

\begin{theorem}
\label{WCP2}
Let $(\A,\Delta)$ be an algebraic quantum group and
 $U$ be a unitary corepresentation of $(\A_r,\Delta_r).$
 The following conditions are equivalent:
\begin{itemize}
\item[(1)] $I \otimes 1 \prec U,$ that is, $U$ has the WCP ;
\item[(2)] there exists $\psi \in S(B(\H_U))$ such that
$\psi((\omega \bar\otimes \iota)U)= \omega (I),
\ \forall \omega \in \M_{*}$;
\item[(3)] there exists a net $(\xi_i)$ of unit vectors in $\H_U$
such that
$$
\lim_i \ \|U(v \otimes v_i) - v \otimes v_i\|_2 = 0, \ \forall v \in \H \ ;
$$
\item[(4)] there exists a net $(\xi_i)$ of unit vectors in $\H_U$
such that
$$
\lim_i \ (U(v \otimes \xi_i),v \otimes \xi_i) = 1,
\ \forall v \in \H, \|v\|_2=1 \ ;
$$
\end{itemize}
Further, any of these conditions implies that $U$ is left-amenable,
right-amenable
and hypertracial.

\end{theorem}

We will illustrate in the next section that
 amenability of $U$ does not imply in general  that $U$ has the  WCP.
 We now collect some elementary facts about the WCP.
 
\begin{proposition}
\label{WCPprop}
Let $(\A,\Delta)$ be an algebraic quantum group and
let $U, V$ be unitary corepresentations of $(\A_r,\Delta_r).$

\vspace{1ex}
\noindent (1) If $U$ has the WCP, then $ \overline{U}$ has also the WCP.

\vspace{1ex}
\noindent
(2) If $U$ and $V$ have the WCP, then $U \times V$ also has the WCP.

\vspace{1ex}
\noindent
(3) If $U$ has the WCP, then $ U_{HS}$ has also the WCP.

\vspace{1ex}
\noindent
(4) If $U \times V$ has the WCP, then $U$ is left-amenable and $V$ is
right-amenable.

\vspace{1ex}
\noindent
(5) If $ U_{HS}$ has the WCP, then
$\overline{U}$ is left-amenable and $U$ is right-amenable.
\end{proposition}

\begin{proof}The proof of assertion (1) is an easy exercise, left to
the reader.
Assertion (2) is a straightforward application of condition (2) in Theorem
\ref{WCP2}.
Assertion (3) follows from (1) and (2).
To prove assertion (4), assume that $U \times V$ has the WCP. The last
assertion in Theorem
\ref{WCP2} tells us then that $U \times V$ is left- and right-amenable.
If now  $M \in S(B(\H_{U})
\otimes B(\H_{V}))$ is a left-invariant (resp. right-invariant) mean
for $U \times V,$ then
one checks without difficulty that $m_{U}(x)= M( x \otimes I_{\H_{V}})$
(resp.
$m'_{V}(y) = M (I_{\H_{U}} \otimes y)$), $x \in B(\H_{U})$
(resp. $y \in B(\H_{V})$), is a left-invariant (resp.
right-invariant) mean for $U$ (resp. $V$). This shows (4). Finally,
assertion (5) follows clearly from (4).

\end{proof}

\emph{Remark.} Let $u$ denote a unitary representation of a discrete
group $\Gamma.$ One of the main results of Bekka in \cite{Bek} is
that $u$ is amenable if, and only if, its associated Hilbert-Schmidt
representation weakly contains the trivial representation. An
interesting question is whether some quantum group version of this result
is true, that is, whether the converse of assertion (5) in
Proposition \ref{WCPprop} holds, at least in some cases. We will
return to this question in Section 9.

\begin{corollary}
Let $U$ be a finite-dimensional unitary corepresentation of $(\A_r,\Delta_r).$
Assume that
$(R \odot  \iota)U = U^*.$ (This is known to hold in the Kac algebra
case, cf.\ \cite[Proposition 1.5.1]{EnS} ).

Then $I \otimes 1 < U_{HS},$
and  $U$ is both left- and right-amenable.
\end{corollary}

\begin{proof}Write
$U = \sum_i a_i \otimes b_i,$
where $a_i \in \M, b_i \in B(\H_U), i = 1 \ldots n$. We use the
notation introduced in Section 4. Using Proposition \ref{UHS}
we may write $U_{HS}=(\iota
\odot l)U (R \odot r)U.$
As $(R \odot r)U = (\iota \odot r)U^*,$
we have $(R \odot r)U = \sum_j a^*_j \otimes r(b^*_j).$
Let $\xi = I_{\H_U} \in HS(\H_U).$
For any    $\eta \in \H$ we get
$$U_{HS}(\eta \otimes \xi)
= ((\iota \odot l)U (R \odot r)U)(\eta \otimes \xi)$$
$$= (\sum_i a_i \otimes l(b_i))(\sum_j a^*_j \otimes r(b^*_j))
(\eta \otimes \xi)$$
$$= (\sum_{ij} a_i a^*_j \otimes l(b_i)r(b^*_j))(\eta \otimes \xi)
= \sum_{ij}  a_i a^*_j \eta \otimes b_i b^*_j$$
$$= (\sum_{ij}  a_i a^*_j \otimes b_i b^*_j) (\eta \otimes \xi)
= \eta \otimes \xi$$
where, in the last equality, we have used the fact that
$\sum_{ij}  a_i a^*_j \otimes b_i b^*_j = UU^* = I_\H \otimes
I_{\H_U}.$ Thus, appealing to Proposition \ref{strongcont}, we have shown
that $I \otimes 1 < U_{HS}.$ We may then apply Proposition \ref{WCPprop}
(5) and conclude that $U$ is right-amenable. Finally, as 
$\overline{U}$ is also finite-dimensional, we then easily deduce that 
$\overline{U}$ is right-amenable, hence that $U$ is left-amenable. 
\end{proof}

\begin{proposition}
\label{Wprop}
Let $(\A,\Delta)$ be an algebraic quantum group and
consider its multiplicative unitary $W$ as a unitary corepresentation of
$(\A_r,\Delta_r).$

\vspace{1ex}
Then $W$ has the WCP if, and only if, $W_{HS}$ has the WCP.
\end{proposition}

\begin{proof}
Using assertion (3) of Proposition \ref{WCPprop}, it suffices
to show that $W$ has the WCP whenever $W_{HS}$ has it. So assume
that $I \otimes 1 \prec W.$
Using the absorbing property of $W$ (cf. our remark after
Proposition \ref{absor}), Ê
we obtain that $W_{HS} =
\overline{ W} \times  W$
is unitarily equivalent to
$I_{\overline{ W}} \times  W
= W_{13}(I_\H \otimes I_{\overline{\H}})_{12}
= W_{13}$.
According to Theorem \ref{WCP2},
there exists $\psi \in S(B(\overline{\H}\otimes\H))$
such that $\psi((\omega \bar\otimes \iota \bar\otimes\iota)W_{13})
= \omega(I), \ \omega \in \hat\M_*$.
Define a state $\psi'$ on  $B(\H)$ by
$\psi'(x) = \psi (I_{\overline\H} \otimes x), \ x \in B(\H).$
Then $$\psi'((\omega \bar\otimes \iota)  W)
= \psi(I_{\overline\H} \otimes (\omega \bar\otimes \iota)  W)
= \psi((\omega \bar\otimes \iota \bar\otimes\iota) W_{13})
= \omega(I).$$ Using Theorem \ref{WCP2} again, we deduce that $W$ has the
WCP,
as desired.
\end{proof}

\begin{corollary}
\label{hatWprop}
Let $(\A,\Delta)$ be an algebraic quantum group and
consider the multiplicative unitary $\hat W$ as a unitary corepresentation
of $(\hat{\A}_r,\hat{\Delta}_{r,\rop}).$
Then

\vspace{1ex}
$(\A,\Delta)$ is co-amenable if, and only if, $\hat{W}_{HS}$ has the WCP.

\end{corollary}

\begin{proof}
 We just have to combine the dual version of
Proposition \ref{Wprop} with Corollary \ref{charcoamen2}.
\end{proof}

Our interest in such a result is that it is presumably easier to
establish that $\hat{W}_{HS}$ has the WCP than to establish that $\hat W$
has the WCP if one wants to show that $(\A,\Delta)$ is co-amenable.

\vspace{2ex}
We conclude this subsection with another proposition involving containment
and amenability.

\begin{proposition} Let $(\A,\Delta)$
denote an algebraic quantum group and $U, V$ be  unitary
corepresentations of $(\A_{r},\Delta_{r}).$

If $U$ is left- (resp.\ right-) amenable and $U < V$, then $V$ is left-
 (resp.\ right-) amenable.
\end{proposition}

\begin{proof}
Let $T \in Mor(U,V)$  be such that $T^*T=I$ and assume that $U$ is
left-amenable.
Define $\tilde{T}: B(\H_V) \to B(\H_U)$ by $\tilde{T}(x)=T^*x T,$
and note that $\tilde{T}$ is a normal $*$-preserving completely positive
unital linear map.
Then note that, since $T \in Mor (U,V)$ and $ T^*T=I, $ we have
\begin{align*}
(\omega \bar\otimes \iota)U
& = T^* (\omega \bar\otimes \iota)V T \\
& = \tilde{T} (\omega \bar\otimes \iota)V \\
& = (\omega \bar\otimes \iota)(\iota\bar\otimes \tilde{T})V
\end{align*}
for all $\omega \in \M_{*}$
(the last equality can be checked for $V$ elementary first and then
by continuity).
As $\{ \ (\omega \bar\otimes \iota) \, | \, \omega \in \M_{*} \, \}$
separate
the elements of
$\M \bar\otimes B(\H_U)$,
we get $U = (\iota\bar\otimes \tilde{T})V$.

\vspace{1ex}
Now, let $m_U \in S(B(\H_U))$ be a left-invariant mean for $U,$
so
$$m_U((\omega \bar\otimes \iota)(U^*(I \otimes y)U))
=\omega(I) m_U(y)$$ for all $y \in B(\H_U)$ and $\omega \in \M_{*}$.

Define $m_V \in S(B(\H_V))$ by $m_V = m_U \tilde{T}$.
Then, for $x \in B(\H_V)$ and $\omega \in \M_{*},$ we get
$$m_V(x) =m_U \tilde{T}(x)
= m_U((\omega \bar\otimes \iota)(U^*(I \otimes \tilde{T}(x))U)).$$
Since $(\iota\bar\otimes\tilde{T})V = U$ is unitary and
$\iota\bar\otimes\tilde{T}$ is completely positive,
it follows from a well known result of M.D. Choi, see e.g. \cite[9.2
]{Str}, that $\iota\bar\otimes\tilde{T}$ is multiplicative at $V$ and
$V^*$. Hence, we get

$$(\iota\bar\otimes\tilde{T})(V^*(I \otimes x)V) =
(\iota\bar\otimes\tilde{T}) V^* (I \otimes \tilde{T}(x))
(\iota\bar\otimes\tilde{T})V .$$
Thus
\begin{align*}
m_V(x)
& = m_U ( (\omega \bar\otimes \iota) (\iota\bar\otimes\tilde{T})
(V^*(I \otimes x) V) ) \\
& = m_U ( \tilde{T} (\omega \bar\otimes \iota) (V^*(I \otimes x) V))
\\
& = m_V( (\omega \bar\otimes \iota ) (V^*(I \otimes x) V) ).
\end{align*}
So $m_V$ is a left-invariant mean for $V$ and $V$ is
left-amenable. The proof of the resp.\ part of the statement is similar.
\end{proof}

It would be interesting to know whether this result still holds if
one replaces strong containment with weak containment.  Bekka has
shown \cite[Corollary 5.3]{Bek} that this is true in the classical case.

\subsection{On property (T)}

We introduce a version of Kazhdan's property (T) \cite{HV} for algebraic
quantum groups.
Then,
as in the classical case,
we show that
every compact quantum group has property (T).
This implies that none of the non-trivial irreducible corepresentations
of a compact quantum group has the
WCP. Furthermore, we show that compactness may be characterized by
having property (T) together with co-amenability of the dual quantum
group.

 \begin{definition}
Let $(\A,\Delta)$ be an algebraic quantum group.
We say that $( \A,  \Delta)$
has property (T) if $I \otimes 1 \prec  U \Rightarrow I \otimes 1 <
U$ for all unitary corepresentations $U$ of
$(\A_{r},\Delta_{r}),$
in other words, if $\hat{\varepsilon}_u \prec \hat{\pi} \Rightarrow
\hat{\varepsilon}_u <
\hat \pi$
for all non-degenerate $*$-representations $ \hat \pi$ of $\hat{\A}_u.$
\end{definition}

\begin{theorem}\label{propT}
Let $(\A,\Delta)$ be an algebraic quantum group of compact type.
Then $( \A,  \Delta)$ has property (T).
\end{theorem}

\begin{proof}
    Let $U$ be a unitary corepresentation of $( \A_{r},
    \Delta_{r})$ and  assume that $U$ has the WCP.
    To show the theorem, we have to show that $\hat{\varepsilon}_u <
    \hat{\pi}_{U}.$

    \smallskip
    Since $( \A,  \Delta)$ is of compact type,  $( \A_{r}, \Delta_{r})$
    is a
    compact quantum group in the sense of Woronowicz. Its Haar state
    $\varphi_{r}$ is then left- and right-invariant, and it has a unique
     extension to a normal state on $\M$ which we also denote by
$\varphi_{r}.$

      \smallskip
Now, let $\xi \in \H_{U}$ and set $\eta
    = ((\varphi_{r}
    \otimes \iota) U) \xi \in \H_{U}.$
    Then, for all $v \in \H$, we have
$$ U ( v \otimes \eta) = U  ( I \otimes ( \varphi_{r}  \otimes
    \iota) U ) (v \otimes \xi),$$
    while
    $$ v \otimes \eta =  (I \otimes ( \varphi_{r}  \otimes
    \iota) U ) ( v \otimes \xi).$$
\noindent But
     $$ I \otimes (( \varphi_{r}  \otimes \iota) U)
    = ( \varphi_{r} (\cdot) I   \otimes  \iota) U$$
$$= (\iota \otimes \varphi_{r} \otimes \iota) ( \Delta_{r}
    \otimes \iota) U \ \ \ \textnormal{ (using invariance of} \
\varphi_{r})$$
$$= (\iota \otimes \varphi_{r} \otimes \iota)
    (U_{13} \, U_{23}) = U (I \otimes ( \varphi_{r}  \otimes \iota)
    U).$$

    Hence, we get $U ( v \otimes \eta) = v \otimes \eta$ for all $v
    \in \H.$

    \smallskip
    Using the dual version of  Proposition \ref{strongcont}, we will then
have shown that
$\hat{\varepsilon}_u < \hat{\pi}_{U}$ if we can show that the vector $\eta$
may be chosen to be
    non-zero. This may be seen as follows.
Since $U$ has the WCP, we  know from Theorem
    \ref{WCP2}
    that there exists a state $\psi$ on $B(\H_{U})$ such that
    $\psi (( \varphi_{r} \bar{\otimes} \iota) U) = 1$. This implies
    that $(\varphi_{r} \bar{\otimes} \iota) U \neq 0.$ Hence, there
    exists at least one $\xi \in \H_{U}$ such that
    $$0 \neq  (( \varphi_{r} \bar{\otimes} \iota) U) \xi = (( \varphi_{r}
   \otimes \iota) U) \xi = \eta,$$
    as desired.
 \end{proof}

 \begin{theorem}
     Let $(\A,\Delta)$ be an algebraic quantum group.
Then $( \A,  \Delta)$ is of compact type if and only if $( \A,
\Delta)$ has property (T) and $(\hat \A, \hat{\Delta})$ is co-amenable.
\end{theorem}

\begin{proof} If $( \A,  \Delta)$ is of compact type, then we know
from Theorem \ref{propT} that  $( \A,
\Delta)$ has property (T). Further, $(\hat \A, \hat{\Delta})$ is then of
discrete type and therefore co-amenable \cite{BMT2}.

Conversely, assume that $( \A,
\Delta)$ has property (T) and $(\hat \A, \hat{\Delta})$ is co-amenable.
Then, using  the dual version of Corollary \ref{charcoamen1}, we get
 $\hat{\varepsilon}_{u} \prec \hat{\pi}_{W}, $ hence
$\hat{\varepsilon}_{u} < \hat{\pi}_{W}. $ This means that there exists
a $T: \Complessi \to \H$ such that $T^{*}T = 1$ and
$\hat{\varepsilon}_{u}(y) = T^{*} \hat{\pi}_{W}(y) T$ for all $y \in
\hat{\A}_{u}.$ Thus we have $\hat{\varepsilon}_{u}(y) = ( \,
\hat{\pi}_{W}(y) \eta, \eta \,)$ for all $y \in
\hat{\A}_{u},$ where $\eta = T(1)$ is a unit vector in $\H.$

Let
$\psi$ denote the vector state $\omega_{\eta}$ on $B(\H).$ Then,
proceeding as
in the proof of Theorem \ref{WCPequiv}, (1) implies (2),  we get
$\psi(\omega \bar \otimes \iota) W = \omega(I_{\H})$ for all $\omega
\in \M_{*}.$ As $\psi$ is normal, this gives $ \omega ( \iota \bar
\otimes \psi) W = \omega(I_{\H})$ for all $\omega
\in \M_{*},$ hence $( \iota \bar
\otimes \psi) W = I_{\H}.$ Since $\A_{r}$ is the norm closure of
$\{ \, ( \iota \bar \otimes \phi) W \, | \,  \phi \in B(\H)_{*} \, \},$ we
get
$I_{H} \in \A_{r},$ that is $(\A_{r}, \Delta_{r})$ is compact, as
desired.
\end{proof}

It is clear that a more detailed study of property (T) for algebraic
quantum groups would be an interesting task (see \cite{PJ} for the
case of Kac algebras). However, we don't
elaborate further on this   as it would take us too
far apart from our main theme in this paper.

\section{Amenability vs. co-amenability vs. WCP}

Let $(\A,\Delta)$
denote an algebraic quantum group and $U$ be a unitary
corepresentation of $(\A_{r},\Delta_{r}).$

\vspace{1ex}
\noindent
We show  that some of the notions introduced  in the previous sections
concerning $U$ are  different from each other by producing counter examples
to
the various possible implications. We consider here only
left-amenability, as we may obtain similar statements for
right-amenability by considering the conjugate of $U.$

\vspace{2ex}
\noindent (1) $U$  co-amenable does not imply that $U$ has the WCP.

\vspace{1ex}
In fact, pick $(\A,\Delta)$  of discrete type and  such that  $(\hat
A, \hat \Delta)$ is not co-amenable (e.g. $\A = C_{c}({\Bbb F}_2)$).
Then $(\A,\Delta)$ is co-amenable since it is of discrete type, cf.
\cite[Theorem 4.1]{BMT2}. Hence, every corepresentation of it is
co-amenable, by
Theorem \ref{allcoamen}. In particular $U$ is co-amenable. But $U$ has not
the WCP since $(\hat
A, \hat \Delta)$ is not co-amenable, by the dual version of Corollary
\ref{charcoamen2}.

\vspace{2ex}
\noindent (2)  $U$ has the WCP does not imply that $U$ is co-amenable.

\vspace{1ex}

Indeed, pick $(\A,\Delta)$ non co-amenable and of compact type (e.\ g.\
$\A = \Complessi[{\Bbb F}_2]).$
Again pick $U = W$. Now,  $(\hat A, \hat{\Delta}_{\rop})$ is co-amenable
(being of
discrete type). Hence, $U$ has the WCP, using the dual version of
Corollary \ref{charcoamen2}.
On the other hand, $U$ is not co-amenable, according to
Theorem \ref{allcoamen}.

\vspace{2ex}
\noindent (3)  $U$ left-amenable does not imply that $U$ is co-amenable.

\vspace{1ex}
Again, pick $(\A,\Delta)$ non co-amenable and of compact type
and let $U=W.$ Since any compact quantum group is amenable,
see the paragraph preceding Theorem 4.7 in \cite{BMT2},
$U$ is left-amenable according to Theorem \ref{allamen}.
On the other hand, according to Theorem \ref{allcoamen}, U is not
co-amenable.

\vspace{2ex}
\noindent (4)  $U$ co-amenable does not imply that $U$ is left-amenable.

\vspace{1ex}
Let $(\A,\Delta)$ be non-amenable and of discrete type.
Being  co-amenable, all its
unitary corepresentations  are then co-amenable.
However, they cannot all be amenable.

 \vspace{2ex}
\noindent (5) $U$  left-amenable does not imply that $U$ has the WCP.

\vspace{1ex}
 Indeed, let $\Gamma$ be  any non-trivial  finite group and let $\A =
C(\Gamma).$
 Let $ \Delta$ be
the usual
 co-product on $\A.$ Then pick a non-trivial
 irreducible unitary representation $u$ of $\Gamma$ and let
 $U$ be the unitary corepresentation of
 $(\A_{r}, \Delta_{r})$ associated with $u.$  Now, it is clear
that  $(\A_{r}, \Delta_{r})$
 is amenable and has property (T) (since it is compact). Then $U$ is
 amenable ( by Theorem \ref{allamen}), but $U$ has not the WCP
 (as remarked at the beginning of subsection 7.3).

 \vspace{2ex}
\emph{Remark.}
Let $(\A, \Delta)$ be of  compact type. As used several
 times by now, $(\A, \Delta)$  is then amenable
and all
the unitary corepresentations of $(\A_r,\Delta_r)$ are therefore amenable.
If $(\A, \Delta)$ is also co-amenable
(e.g. we may take the compact matrix pseudogroup $\A=$
SU$_q(2),$ cf.
\cite{Ba2, BMT1}),
all these corepresentations are then also co-amenable.
Further, as $(\A, \Delta)$ has  property T,
 we get that
none of the non-trivial irreducible corepresentations of $(\A_r,\Delta_r)$
satisfies the WCP.

 On the other hand,
$(\hat\A,\hat\Delta)$ is always co-amenable
since it is of discrete type.
Hence,  all the unitary corepresentations of
$(\hat\A_r,\hat\Delta_{r,\rop})$
are co-amenable. If $(\A, \Delta)$ is also co-amenable, then we know
that $(\hat \A, \hat \Delta)$  is amenable, hence all these
corepresentations are then also amenable.

\section{Amenability and discrete quantum groups }

As we pointed out in connection with Proposition \ref{WCPprop}, it would be
interesting to know whether the converse of Proposition
\ref{WCPprop} (5) holds, that is, whether the right-amenability of
$U$ implies that $U_{HS}$ has the WCP.  It does in the classical case,
and this is one of the
major result in \cite{Bek}. The problem of going from amenability of $U$
to the WCP for $U_{HS}$ seems much more delicate in the general case.

We will now present a proof which works for an algebraic quantum group
of discrete type having  a (compact) dual with a tracial Haar state.
As a consequence, we obtain a new proof of the fact that amenability
of such a quantum group is equivalent to the co-amenability of its
dual, which has been previously established by  Ruan \cite[Theorem
4.5]{Ruan},
see also \cite{BMT3}.

\vspace{2ex}
For notational reasons, we let $(\A, \Delta)$ be an algebraic quantum
group of compact type and consider its dual $(\hat \A, \hat \Delta)$
which is then of discrete type. We use the description of
$(\hat \A, \hat \Delta)$ given in Proposition \ref{discrete} and the
notation introduced there. We
denote by $\tS = \hat{S}_{\rop}$ the antipode of $(\hat A,
\hat{\Delta}_{\rop}),$ and by
$\tR$ the anti-unitary  antipode of $(\hat{\A}_{r},
\hat{\Delta}_{r,\rop})$ (which is defined on $\hat\M$). For each $\a \in
A,$ we
denote  the  central minimal projection of $\hat\M$ which
is given by $\hat\pi_{r}(p_{\a})$ with the same symbol $p_{\a}$.
Further, we identify $\hat\pi_{r}( \hat{\A}_{\a}) = p_{\a} \hat\pi_{r}(
\hat{\A})$ with $ \hat{\A}_{\a} = M_{d_{\a}}(\Complessi)$ and
let $Tr_{\a}$ denote its canonical trace. Finally, we denote
    the canonical injection from $\hat\A$ into $\H$ by  $\hat\Lambda.$

\vspace{1ex}
Let now $U$ be a unitary corepresentation of $(\hat{A}_{r},
\hat{\Delta}_{r,\rop}).$
We remark  that, using the
above identifications and the properties of $p_{\a}$, one easily
deduces
that
 $$(p_{\a} \otimes I)U = U(p_{\a} \otimes I) \in \hat{\A}_{\a} \odot
 B(\H_{U}) , \ \ U(p_\alpha \otimes y)U^{*} \in \hat{\A}_{\a} \odot
 HS(\H_{U})$$
for all $\a \in A$ and $y \in HS(\H_{U}).$
We denote by $T_{\a}$ the trace on  $\hat{\A}_{\a} \odot
 B(\H_{U})$ given by $T_{\a} = Tr_{\a} \odot Tr,$ where $Tr$
denotes the canonical trace on $HS(\H_{U}).$ Further, we denote by
$\| \cdot \|_{1,\a}$ and $\| \cdot \|_{2,\a}$ the associated norms on
$\hat{\A}_{\a} \odot TC(\H_{U})$ and $\hat{\A}_{\a} \odot
 HS(\H_{U}),$ respectively.

 \vspace{1ex}
We establish  a series of lemmas.

\begin{lemma} \label{WCPdisc}
    For each $\a \in A,$ set
    $$b_{\a} = M_{\alpha} \sum_{i,j=1}^{d_{\alpha}} f_{1}(u^{\alpha}_{ji})
    u^{\alpha}_{ij}, $$
    so that we have $\hat{b}_{\a} = p_{\a}.$
The following conditions are equivalent:

   \vspace{1ex}
    \noindent (1) $U$ has the WCP.

    \vspace{1ex}
    \noindent (2) There exists a net $(\xi_{i})$ of unit vectors in
    $\H_{U}$ such that
    $$\lim_{i} \| U(\hat\Lambda(p_{a}) \otimes \xi_{i}) -
    (\hat\Lambda(p_{a}) \otimes \xi_{i})\|_{2} = 0 \ \  \forall \a
    \in A.$$

    \vspace{1ex}
    \noindent (3) There exists a state $\phi$ on $B(\H_{U})$ such that
    $$ \phi (( \omega_{\hat{\Lambda}(p_{\a})} \bar \otimes \iota) U) =
    \omega_{\hat{\Lambda}(p_{\a})} (I), \ \ \forall \a \in A.$$

    \vspace{1ex}
    \noindent (4) There exists a state $\phi$ on $B(\H_{U})$ such that
    $$ \phi \pi_{U}(b_{\a}) = M_{\a}^2 , \ \ \forall \a \in A.$$

    \vspace{1ex}
    \noindent (5) There exists a state $\phi$ on $B(\H_{U})$ such that
    $ \phi \pi_{U} = \e_{u}.$

\end{lemma}

\begin{proof} (1) $\Rightarrow$ (2) and (2) $\Rightarrow$ (3) follow
as in the proof of Theorem \ref{WCPequiv}.

\vspace{1ex}
\noindent (3) $\Rightarrow$ (4) : Assume that (3) holds and let $\phi$ be as
in (3).
We will show that $\phi$ satisfies (4). Fix $\a \in A.$

We first observe that $S(b_{\a}^*) = b_{\a}.$ Indeed,
$$S(b_{\a}^*) = M_{\a} \sum_{i,j} \overline{f_{1}(\u_{ji})}
S((\u_{ij})^*)$$
$$= M_{\a} \sum_{i,j} \overline{f_{1}(\u_{ji})}
\sum_{k,l}f_{1}(\u_{jk})f_{-1}(\u_{li}) \u_{kl}$$
$$= M_{\a} \sum_{i,j,k,l} f_{-1}((\u_{ji})^*) f_{-1}(\u_{li})
f_{1}(\u_{jk}) \u_{kl}$$
$$= M_{\a} \sum_{j,k,l} f_{-1}(\sum_{i}(\u_{ji})^*\u_{li})
f_{1}(\u_{jk}) \u_{kl}$$
$$= M_{\a} \sum_{j,k,l} f_{-1}(\delta_{jl}I)
f_{1}(\u_{jk}) \u_{kl} = M_{\a} \sum_{k,l} f_{1}(\u_{lk}) \u_{kl} =
b_{\a}.$$
Thus, we have
$$\hat{b}_{\a} (S((b_{\a})^*))^{\wedge} = p_{\a} p_{\a} = p_{\a}=
(S((b_{\a})^*))^{\wedge}.$$
Therefore, using the result from \cite{BMT3} recalled at the beginning
of subsection 3.1, we have
$$Q(b_{\a}) = \omega_{\Lambda(b_{\a}),\Lambda(b_{\a})} =
\omega_{\hat{\Lambda}(p_{\a})}.$$
Using \cite[Theorem 3.2]{BMT3} and the assumption that $\phi$ satisfies (3),
we get
$$ \phi (\pi_{U} ( b_{\a}))= \phi ( (Q(b_{\a}) \bar\otimes \iota)U)$$
$$= \phi ( ( \omega_{\hat{\Lambda}(p_{\a})} \bar\otimes \iota)U)
= \omega_{\hat{\Lambda}(p_{\a})} (I) = \hat\psi(p_{\a}^*
p_{\a})$$
$$= \hat\psi(p_{\a}) = \hat\psi(\hat{b}_{\a}) = \e(b_{\a})$$
$$= M_{\alpha} \sum_{i,j} f_{1}(u^{\alpha}_{ji})
    \e(u^{\alpha}_{ij}) = M_{\a} \sum_{i}f_{1}(u^{\alpha}_{ii}) =
    M^2_{\a},$$
    and (4) is proved.

    \vspace{1ex} \noindent
    (4) $\Rightarrow$ (5) : Assume (4) holds and let $\phi$ be as in (4).
    Let $\eta_{u}$ be the state on $\A_{u}$
   given by $\eta_{u} = \phi \pi_{U}$ and let $\eta$ denote the
   restriction of  $\eta_{u}$ to $\A.$ To show that $\phi$ satisfies
   (5), that is $\eta_{u} =\e_{u},$ it suffices to show that $\eta = \e.$

   \vspace{1ex}
   Fix $\alpha \in A.$
   As $\sum_{i,j} f_{1}(u^{\alpha}_{ji}) u^{\alpha}_{ij} = \sum_{i}
   f_{1}*\u_{ii},$ we have $$(*) \ \ \ \ \phi(\sum_{i} f_{1}*\u_{ii})
=
   M_{\a}= f_{1 }(\sum_{i}\u_{ii}).$$

   Set $d= \sum_{i}\u_{ii}$ and observe that we may write $(*)$ as
$\eta
   f_{1} (d) = f_{1}(d).$

   \vspace{1ex}
   Now, set $ X_{ij}= f_{1/2}*\u_{ij} - f_{1/2}(\u_{ij})I \, \in \A.$ Then

   \vspace{1ex}
   $\eta(\sum_{i,j}X^*_{ij}X_{ij})$
   $$= \sum_{i,j} \eta((f_{1/2}*\u_{ij})^*(f_{1/2}*\u_{ij})) - 2 Re(
   \overline{f_{1/2}(\u_{ij})} \eta(f_{1/2}*\u_{ij})) +
   |f_{1/2}(\u_{ij})|^2$$
   $$=\sum_{i,j}\eta(\sum_{k,l}(\u_{ik})^*
\overline{f_{1/2}(\u_{kj})} \u_{il}f_{1/2}(\u_{lj})$$
$$- 2 Re (
\overline{f_{1/2}(\u_{ij})}
\sum_{k}f_{1/2}(\u_{kj}) \eta(u_{ik}))  + |f_{1/2}(\u_{ij})|^2$$
$$ =\sum_{j,k,l} (\delta_{kl}\eta(I) \overline{f_{1/2}(\u_{kj})}
f_{1/2}(u_{lj})) - 2 Re ( \sum_{i,j,k} \eta(u_{ik})f_{1/2}(\u_{ji})
f_{1/2}(\u_{kj}) ) + \sum_{i}f_{1}(\u_{ii})$$
$$=\sum_{j,k} f_{1/2}(\u_{jk})f_{1/2}(\u_{kj}) - 2 Re (\eta
f_{1}(\u_{ii})) + \sum_{i}f_{1}(\u_{ii})$$
$$= 2 \sum_{i}f_{1}(\u_{ii}) - 2 Re (\eta
f_{1}(\u_{ii})) = 2 (f_{1}(d) - Re(\eta f_{1} (d))).$$
Now, as $\eta f_{1}(d) = f_{1}(d) = M_{\a}$ is  real, we
get $$\eta(\sum_{i,j}X^*_{ij}X_{ij})= 0.$$
Since $\eta$ is a positive linear functional, this implies that
$\eta(X^*_{ij}X_{ij}) = 0$ for all $i, j.$ Using the Cauchy-Schwarz
inequality, we obtain  $\eta(X_{ij}) = 0$ for all $i, j,$ that is,
$$\eta (f_{1/2}* \u_{ij}) = f_{1/2}(\u_{ij}) , \ \ \forall i, j,$$
hence
$$\eta (f_{1/2}* a) = f_{1/2}(a) , \ \ \forall a \in \A,$$
by linearity.

For any $b \in \A,$ we let $a=f_{-1/2}*b$ and apply the above. this
gives
$$ \eta(f_{1/2} * f_{1/2} * b) = f_{1/2}(f_{-1/2}*b),$$
that is, $\eta(b) = \e(b).$ Thus, we have shown that $\eta = \e,$ as
desired.

\vspace{1ex}
\noindent (5) $\Rightarrow$ (1) : Assume (5) holds. Then we clearly
have Ker $\pi_{U} \subset $ Ker $\e_{u},$ that is, $U$ has the WCP.

    \end{proof}

\begin{lemma}\label{SU}
    Let $\a \in A$ and set $p_{\beta} = \tS(p_{\a}).$ Then we have
    $$ (\tS \odot \iota) (( p_{\a} \otimes I) U) = U^* (p_{\beta}
    \otimes I).$$
    \end{lemma}

    \begin{proof} From the proof of \cite[Proposition 3.4]{KV}, we
    know that $$\tS (( \iota \bar \otimes \omega)U) = (\iota \otimes
    \omega) (U^{*}), \ \omega \in B(\H_{U})_{*}.$$
    Hence, we get
    $$( \iota \bar \otimes \omega)(\tS \odot \iota) (( p_{\a} \otimes I)
    U) = \tS ( p_{\a}( \iota \bar \otimes \omega)U)) $$
    $$= \tS (( \iota \bar \otimes \omega)U) p_{\beta} = (\iota \otimes
    \omega) (U^{*})p_{\beta} $$
    $$ = ( \iota \bar \otimes \omega) (U^* (p_{\beta} \otimes I))$$
    for all $\omega \in B(\H_{U})_{*}.$
    \end{proof}

\begin{lemma}\label{Bek1}
Assume that $U$ is right-amenable.
Then there exists a net $(y_i)$  in $ \{ y \in HS(\H_U) \, | \,
y \geq 0, \, \|y\|_2 = 1 \}$
such that
$$\lim_i
\|U(p_\alpha \otimes y_i)U^{*} - p_\alpha \otimes y_i\|_{2,\a}
= 0, \ \forall \alpha \in A.$$
\end{lemma}

\begin{proof} We begin as in the proof of Proposition \ref{BekTC}.
 Let $m'_{U}$ be
 a right-invariant mean for $U,$ so
 $$m'_U((\omega \bar\otimes \iota)U(I \otimes x)U^*) =
 \omega(I)m'_U(x)$$
 for all $x \in B(\H_{U}), \omega \in \hat{\M}_{*}.$
 As the normal states are
 weak$^{*}$-dense in $S(B(\H_{U}))$, we may pick a net $(s_{i}) \subset
 TC(\H_{U})^{+}_{1}$ such that $m'_{U}$ is  a weak$^{*}$-limit point of
 the net $(Tr(s_{i} \, \cdot)) \subset S(B(\H_{U})).$

 Now, we define
 a net $(y_i)$  in $ \{ y \in HS(\H_U) \, | \,
y \geq 0, \, \|y\|_2 = 1 \}$ by setting $y_{i}= s_{i}^{1/2}$ for all i.

\vspace{1ex}
 Let $\a \in A.$ Hereafter, we write  $\hat{a}_\alpha$ to
 denote $\hat{a}p_\alpha \in \hat{\A}_{\a}$ whenever $ a \in \A.$

 Let $b, b' \in \A.$ Set $c_\alpha =
\hat{\rho}^{-1}(\hat{b'}_\alpha)\hat{b}^*_\alpha.$
Then we have
$$\lim_i
\, \omega_{\hat\Lambda(\hat{b'}_\alpha),\hat\Lambda(\hat{b}_\alpha)}(I)
\, Tr(xs_i)
$$
$$= \lim_i
Tr((
\omega_{\hat\Lambda(\hat{b'}_\alpha),\hat\Lambda(\hat{b}_\alpha)}\bar\otimes
\iota)
(U(I \otimes x)U^*)s_i)$$
$$= \lim_{i} \,(\hat\psi \odot Tr)((\hat{b}^*_\alpha \otimes y_i)U(I \otimes
x)U^*
(\hat{b'}_\alpha \otimes y_i)).$$
Thus
$$\lim_i (\hat\psi \odot Tr)
(\hat{b}^*_\alpha \hat{b'}_\alpha \otimes x y^2_i -
(\hat{b}^*_\alpha \otimes y_i) U(I \otimes x)U^* (\hat{b'}_\alpha \otimes
y_i)) = 0,$$
which gives
$$\lim_i (\hat\psi \odot Tr)
(c_{\a} \otimes x y^2_i -
(c_{\a} \otimes y_i^2) U(I \otimes x)U^*) = 0.$$

Now, write $(p_{\a} \otimes I)U =\sum_r a_r \otimes x_r \in
\hat{\A}_{\a} \odot B(\H_{U}).$

\vspace{1ex}
Then, using that $Tr$ is a trace at the third step and Proposition \ref{discrete}
at the final step, we get that
$$(\hat\psi \odot Tr)(c_\alpha \otimes x y^2_i
- (c_\alpha \otimes  y^2_i)U(I \otimes x)U^*)$$
$$= (\hat\psi \odot Tr)(c_\alpha \otimes x y^2_i )
- \sum_{r,s}(\hat\psi \odot Tr) (c_\alpha a_r a^*_s \otimes y^2_i x_r x x^*_s)$$
$$= (\hat\psi \odot Tr)(c_\alpha \otimes x y^2_i )
- \sum_{r,s} (\hat\psi \odot Tr)( c_\alpha a_r a^*_s \otimes x x^*_s y^2_i x_r)$$
$$=(\hat\psi \odot Tr) ((c_\alpha \otimes x)(p_\alpha \otimes y^2_i
-  \sum_{r,s} (a_r a^*_s \otimes x^*_s y^2_i x_r))$$
$$= (\hat\psi \odot Tr)
((c_\alpha \otimes x)(p_\alpha \otimes y^2_i
-  \sum_{r,s} a_r p_\alpha a^*_s \otimes x^*_s y^2_i x_r)) $$
$$= (Tr_{\a} \odot Tr)
((c_\alpha \otimes x)(p_\alpha f_{-1}\otimes y^2_i
-  \sum_{r,s} a_r p_\alpha a^*_s f_{-1}\otimes x^*_s y^2_i x_r)) $$
 converges to zero. Note that if we let  $b$ and $ b'$  vary in $\A,$
then $c_{\a}$ will give all elements in $\hat{\A}_{\a}$ (using that
$\hat{\A}^2 = \hat{\A}$ and Proposition \ref{discrete}). Hence,  
adapting Namioka's argument \cite[Proof of Theorem 2.4.2]{Gr} by 
considering the locally convex product space $\prod \{  
\hat{\A}_{\a}\otimes TC(\H_{U}) , \a \in A \, \}$ with the product of 
the $\| \cdot \|_{1,\a}$-norm topologies,
we may in fact assume that
 $$\lim_{i} \, \| p_\alpha f_{-1}\otimes y^2_i
- \sum_{r,s} a_r p_\alpha a^*_s f_{-1}\otimes x^*_s y^2_i x_r
\|_{1,\a} = 0.$$

Now, as $\hat{\A}_{\a}$ is a matrix algebra, we can apply any
linear map on the first tensor factor
and still keep convergence in 1-norm.
Doing this with $\tS(\cdot \, f_{1}) \odot \iota,$
we get
$$\lim_{i} \, \|p_\beta \otimes y^2_i
- \sum_{r,s} \tS(a^*_s) p_\beta \tS(a_r)\otimes x^*_s y^2_i x_r
\|_{1,\a} =0.$$

Now, Lemma \ref{SU} says that $ (\tS \odot \iota) ( p_{\a} \otimes I) U) =
U^* (p_{\beta}
    \otimes I).$ Using this,
our last equation reads as
$$ \lim_{i} \|p_\beta \otimes y^2_i - U(p_\beta \otimes y^2_i)U^*\|_{1,\a} =
0.$$
Using the Powers-St\o rmer inequality (see \cite[Lemma 4.2]{Bek}), one gets
then
$$ \lim_{i} \|p_\beta \otimes y_i - U(p_\beta \otimes y_i)U^*\|_{2,\a}= 0.$$
As $\a \to \beta$ is a bijection of $A,$ we are done.
\end{proof}

\begin{lemma} \label{Bek2}
     Assume that $U$ is right-amenable and
     that $(\A,\Delta)$ has a tracial Haar state.
Then there exists a net $(y_i)$  in $ \{ y \in HS(\H_U) \, | \,
y \geq 0, \, \|y\|_2 = 1 \}$
such that
$$\lim_i
\|U_{HS}(\hat\Lambda(p_\alpha) \otimes y_i) - \hat\Lambda(p_\alpha) \otimes
y_i\|_{2}
= 0 \ \ \forall \alpha \in A.$$
\end{lemma}

\begin{proof}Since we assume that $(\A,\Delta)$ has a tracial Haar
state, it is well known that
    $(\M, \Delta_{r})$ is a compact Kac algebra. Hence, $(\hat \M,
    \hat{\Delta}_{r,\rop})$ is then  a discrete Kac algebra and
    we may  identify $\tR$ with the extension of $\tS$ to $\hat \M.$
Recall from Proposition
\ref{UHS} that
$$U_{HS}=(\iota \bar\otimes l)U (\tR \bar\otimes r)U.$$
Let $\a \in A,  \, y \in HS(\H_{U})$ and set $p_{\beta}= \tR(p_{\a}).$ 
Then 
 $$U_{HS}(\hat\Lambda(p_{\a}) \otimes y) =
 (\iota \bar\otimes l)U (\tR \bar\otimes r)(U) (p_{\a}\otimes I) 
 (\hat\Lambda(p_{\a}) \otimes y). $$
Now, using Lemma \ref{SU}, we get
 $$ (\tR \bar\otimes r)(U) (p_{\a} \otimes I) =
 (\tR \bar\otimes r)((p_{\beta} \otimes I)U)$$
 $$= (\iota \odot r)(\tR 
 \odot \iota)((p_{\beta} \otimes I)U) = (\iota \odot r) (U^*(p_{\a} 
 \otimes I)).$$ 
Hence, it follows that
$$U_{HS}(\hat\Lambda(p_{\a}) \otimes y) = 
(\iota \bar\otimes l)U (\iota \odot r)(U^* (p_{\a}\otimes I)) 
 (\hat\Lambda(p_{\a}) \otimes y)$$
$$=(\iota \odot l)(U(p_{\a}\otimes I)) (\iota \odot r)(U^* (p_{\a}\otimes I)) 
 (\hat\Lambda(p_{\a}) \otimes y)=(\hat\Lambda \odot \iota)(U(p_{a}\otimes y)U^*).$$

Now, since $\f$ is assumed to be tracial, we have $f_{1} = \e.$
According to Proposition \ref{discrete}, we then have
$$ \hat \psi (x) = \oplus_{\a} Tr_{\a}(p_{\a}x), \ x \in \hat \A.$$
 It follows that the Hilbert space norm on $\H \otimes HS(\H_{U})$ agrees on
 each subspace $\hat\Lambda(\hat{\A}_{\a}) \odot HS(\H_{U})$ with the $\|
\cdot \|_{2,\a}$-norm  on $\hat{\A}_{\a} \odot HS(\H_{U}).$ Therefore,
 choosing the net $(y_{i})$ as the one provided by Lemma \ref{Bek1}, we get
 $$\lim_i
\|U_{HS}(\hat\Lambda(p_\alpha) \otimes y_i) - \hat\Lambda(p_\alpha) \otimes
y_i\|_{2}
=\|(\hat\Lambda \odot \iota )(U(p_{a}\otimes y_{i})U^*) -
(\hat\Lambda \odot \iota)(p_\alpha \otimes y_i)\|_{2}$$
$$=\lim_{i} \|U(p_\alpha \otimes y_i)U^{*} - p_\alpha \otimes y_i\|_{2,\a}
= 0 $$
for all $\alpha \in A,$ which shows the lemma.
\end{proof}

We are now in position to derive the following analog of \cite[Theorem 5.1]{Bek}.

\begin{theorem}\label{Bek3}
    Assume that $(\A,\Delta)$ is of compact type and has a tracial Haar
state.
    Let  $U$ be a unitary corepresentation of $(\hat A_{r},
\hat{\Delta}_{r,\rop})$.
    Then $U$ is right-amenable if, and only if, $U_{HS}$ has the WCP.
    \end{theorem}

    \begin{proof}
 Assume that $U$ is right-amenable. Combining Lemma \ref{Bek2} with
 Lemma \ref{WCPdisc}, we deduce that $U_{HS}$ has the WCP. The
 converse implication is shown in Proposition \ref{WCPprop} (5).
 \end{proof}
 
 \noindent \emph{Remark}. Let $(\A,\Delta)$ and $U$ be as in Theorem \ref{Bek3}. 
 Recall from 4.2 that we can  associate  with $U$ another 
 Hilbert-Schmidt corepresentation $U_{HS'} \simeq (\overline{U})_{HS}.$   
 As $U$ is left-amenable if, and only if, $\overline{U}$ is 
 right-amenable, we deduce from Theorem \ref{Bek3}  that 
  $U$ is left-amenable if, and only if, $U_{HS'}$ has the WCP.   

\vspace{1ex}
  As an application of Theorem \ref{Bek3} we give  a new proof of the following
 result  (see \cite[Theorem 4.5]{Ruan}, \cite[Theorem 1.1]{BMT3}).

 \begin{corollary} \label{Ruan}
     Assume that $(\A,\Delta)$ is of compact type and has a tracial Haar
state.
     Then $(\A,\Delta)$ is co-amenable if, and only if, $(\hat \A ,
     \hat \Delta)$ is amenable.
     \end{corollary}

     \begin{proof} We know from \cite[Theorem 4.7]{BMT3} that $(\hat \A ,
     \hat \Delta)$ is amenable whenever $(\A,\Delta)$ is co-amenable.
     Assume now that  $(\hat \A, \hat \Delta)$ is amenable. From
     Theorem \ref{allamen}, we deduce that $\hat W$ is right-amenable.
    Using Theorem \ref{Bek3}, we obtain that $\hat W_{HS}$ has
  the WCP.  It follows  from Corollary \ref{hatWprop}
 that $(\A,\Delta)$ is co-amenable.
 \end{proof}

 The question whether the traciality assumption in Corollary
 \ref{Ruan} may be removed remains elusive. It relies on whether the
 traciality assumption in Lemma \ref{Bek2} may be removed.

\bigskip \vskip1cm

\noindent{\bf Acknowledgements.}
Part of this work has been made while R.\ C.\ was visiting the
Department of Mathematics at the University of Oslo in May 2000 and June
2001 and the Centre for Advanced Studies in Oslo in October 2001.
He thanks the members of the Operator Algebras team there
for their hospitality.

\bigskip \smallskip
{\parindent=0pt Addresses of the authors:

\smallskip Erik B\'edos, Institute of Mathematics, University of
Oslo, \\
P.B. 1053 Blindern, 0316 Oslo, Norway.\\ E-mail: bedos@math.uio.no. \\

\smallskip \noindent
Roberto Conti, Mathematisches Institut,\\ Friedrich-Alexander Universit\"at
Erlangen-N\"urnberg\\
Bismarckstr. 1 1/2,  D-91054 Erlangen, Germany. \\ E-mail:
conti@mi.uni-erlangen.de. \\

\smallskip \noindent
Lars Tuset, Faculty of Engineering, Oslo University College, \\
Cort Adelers Gate 30, 0254 Oslo, Norway. \\ E-mail: Lars.Tuset@iu.hio.no.
}


\bigskip

\begin{thebibliography}{42}

\bibitem{BS} Baaj, S. and Skandalis, G.: Unitaires multiplicatifs et
 dualit\'e pour les produits crois\'es de C*-alg\`ebres.
Ann. Sci. Ec. Norm. Sup. 26 (1993), 425--488.

\bibitem{Ba} Banica, T.: Representations of compact quantum groups
 and subfactors. J. Reine Angew. Math. 509 (1999), 167--198.

\bibitem{Ba2} Banica, T.: Fusion rules for representations of
 compact quantum groups. Expo. Math. 17 (1999), 313--338.

\bibitem{Bed} B\'edos, E.: Notes on hypertraces and C*-algebras. J. Operator
Th. 34 (1995), 285--306.

\bibitem{BMT1} B\'edos, E., Murphy, G.J. and Tuset, L.: Co-amenability of
compact quantum groups. J. of Geom. and Phys. 40 (2001), 130-153.

\bibitem{BMT2}    B\'edos, E.,  Murphy ,G. J. and Tuset, L.:
Amenability and co-amenability for algebraic quantum groups.
Preprint (2001).

\bibitem{BMT3}   B\'edos, E.,  Murphy ,G. J. and Tuset, L.:
Amenability and co-amenability for algebraic quantum groups, II.
Preprint(2001).

\bibitem{Bek} Bekka,  M. E. B.:
Amenable unitary representation of locally compact groups,
Invent. Math. 100 (1990), 383--401 .

 \bibitem{Bl} Blanchard, E.: D\'eformations de C*-Alg\`ebres de Hopf.
Bull. Soc. Math. France 124 (1996), 141--215.

\bibitem{Dix} Dixmier,  J.: Les C*-alg\`ebres et leurs repr\'esentations,
Gauthiers-Villars, Paris (1969).

\bibitem{ER} Effros, E.G. and Ruan, Z.-J.: Discrete quantum groups. I. 
The Haar measure. Intern. J. Math. 5 (1994), 681-723.

\bibitem{EnS}  Enock, M. and Schwartz, J.M. : Kac algebras and duality 
of locally compact groups. Springer Verlag, 1992.

\bibitem{ES} Enock, M. and Schwartz, J.M. : Algebres de Kac
moyennables. Pac. J. Math. 125 (1986), 363--379.

\bibitem{Fell} Fell, J.M.G. :
Weak containment and induced representations of groups.
Can. J. Math. 14 (1962), 237--268.

\bibitem{GLR} Ghez, Lima, R., Roberts, J. E. : $W^{*}$-categories. Pac.
J. Math. 120 (1985), 79--109.

\bibitem{Gr} Greenleaf, F. P.: Invariant means on topological groups. 
Van Nostrand (New York), 1969.

\bibitem{HV} de la Harpe, P. and Valette, A.: La propri\'et\'e (T) de
Kazhdan pour
les groupes localement compacts. Asterisque 175, Soc. Math. de
France, 1989.

\bibitem{Ku} Kustermans, J.: Examining the dual of an algebraic quantum
group. Preprint Odense
Universitet (1997). (arXiv:funct-an/9704006).


\bibitem{Kus} Kustermans, J.: Universal C*-algebraic  quantum
groups arising from algebraic quantum groups. Preprint Odense
Universitet (1997). (arXiv:funct-an/9704004).

\bibitem{KuVD} Kustermans, J. and Van Daele, A.: C$^{*}$-algebraic quantum
groups arising from algebraic quantum groups. Int. J. Math. 8
(1997), 1067--1139 .

\bibitem{KV} Kustermans, J. and Vaes, S.: Locally compact quantum
qroups.  Ann. Scient. Ec. Norm. Sup. 33 (2000), 837--934.

 
\bibitem{LR} Longo, R. and Roberts, J. E. : A theory of dimension.
$K$--Theory 11 (1997), 103--159.

\bibitem{Ng1} Ng, C.- K. : Amenability of Hopf C*-algebras.
Proceedings of the 17th OT Conference (2000), 269--284.

\bibitem{Ng2} Ng, C.- K. : Examples of amenable Kac algebras.
Preprint (2000).
 
\bibitem{Pa} Paterson, A. L.: Amenability. Math. Surveys and Monographs
29, Amer. Math. Soc. (1988).

\bibitem{PJ} Petrescu, S. and Joita M.: Property (T) for Kac algebras.
Rev. Roum. Math. Pures Appl. 37 (1992), 163--178.

\bibitem{RoTu} Roberts, J. E. and Tuset, L.:
On the equality of $q$-dimension and intrinsic dimension,
 J. Pure Appl. Alg.  156 (2001), 329--343 .

\bibitem{Ruan} Ruan, Z.-J. : Amenability of Hopf von Neumann algebras
and Kac algebras. J. Funct. Anal. 139 (1996), 466--499.

\bibitem{Str} Stratila, S. : Modular theory in operator algebras.
Abacus Press, Tunbridge Wells, Kent (1981).
 
\bibitem{VD1} Van Daele, A.: Multiplier Hopf algebras. Trans. Amer.
Math. Soc. 342 (1994), 917--932.

\bibitem{VD2} Van Daele, A.: An algebraic framework for group duality.
 Adv. Math. 140 (1998), 323--366.

\bibitem{Vo} Voiculescu, D.: Amenability and Katz Algebras.
Alg\`ebres d'Op\'erateurs et leurs Applications en Physique
Math\'ematique, N 274, Colloques Internationaux C.N.R.S. (1977), 451--457.

\bibitem{Wo1} Woronowicz, S.L.: Compact matrix pseudogroups.
Comm. Math. Phys. 111 (1987), 613--665.

\bibitem{Wo2} Woronowicz, S.L.: Compact quantum groups, in
``Sym\'etries Quantiques'', North Holland, Amsterdam (1998),
845--884.



\end{thebibliography}
\end{document}